\documentclass[preprint,3p,12pt]{elsarticle}
\usepackage{amsmath}
\usepackage{amssymb}
\usepackage{booktabs}
\usepackage{multirow}
\usepackage{float}
\usepackage{subfigure}
\usepackage{lineno,hyperref}
\usepackage{geometry}
\usepackage{color}
\modulolinenumbers[5]

\begin{document}
\begin{frontmatter}
\title{A novel method for constructing high accurate and robust WENO-Z type scheme}
\author[1,2]{Yiqing Shen}
\ead{yqshen@imech.ac.cn}

\author[1,2]{Ke Zhang}

\author[1,2]{Shiyao Li}

\author[1,2]{Jun Peng\corref{cor}}
\ead{pengjun@imech.ac.cn}

\cortext[cor]{Corresponding author}

\address[1]{State Key Laboratory of High Temperature Gas Dynamics, Institute of Mechanics, Chinese Academy of Sciences, Beijing 100190, China}
\address[2]{School of Engineering Science, University of Chinese Academy of Sciences, Beijing 100049, China}

\begin{abstract}

A novel method for constructing robust and high-order accurate weighted essentially non-oscillatory (WENO) scheme is proposed in this paper.  The method is mainly based on the WENO-Z type scheme, in which, an eighth-order global smoothness indicator (the square of the approximation of the fourth-order derivative on the five-point stencil used by the fifth-order WENO scheme) is used, and in order to keep the ENO property and robustness, the constant 1 used to calculate the un-normalized weights is replaced by a function of local smoothness indicators of candidate sub-stencils. This function is designed to have following adaptive property: if the five-point stencil contains a discontinuity, then the function approaches to a small value, otherwise, it approaches to a large value. Analysis and numerical results show that the resulted WENO-Z type (WENO-ZN) scheme is robust for capturing shock waves and, in smooth regions, achieves fifth-order accuracy at first-order critical point and fourth-order accuracy at second-order critical point.

\end{abstract}


\end{frontmatter}

\section{Introduction}\label{Introduction}

Weighted essentially non-oscillatory (WENO) finite difference schemes have been widely studied and applied in computational fluid dynamics. The first WENO scheme was proposed by Liu et al. \cite{Xudong1994} in 1994. Its basic idea is to use a weighted convex combination of fluxes on all candidate sub-stencils instead of the one on the smoothest sub-stencil in ENO scheme \cite{Harten1983}. By assigning adaptive weight to each sub-stencil, the WENO scheme can achieve high order accuracy in smooth regions while keeping ENO property near discontinuities. In \cite{Jiang1996}, Jiang and Shu introduced a general method for calculating smoothness indicators of stencils. Then, Balsara and Shu \cite{Balsara2000} extended the WENO schemes up to 11th order of accuracy. Gerolymos et al. \cite{Gerolymos2009} further developed very-high-order WENO schemes.

In \cite{Henrick2005}, Henrick et al. derived the necessary and sufficient conditions on the weights for fifth-order convergence of a fifth-order WENO scheme and revealed the WENO implemented by Jiang and Shu (WENO-JS) is only third-order accurate at critical points. Then, they proposed a mapping function to correct the weights, which satisfy the sufficient condition for fifth-order convergence even at critical points, resulting in the WENO-M scheme. Borges et al. \cite{Borges2008} introduced a global smoothness indictor (GSI) of higher order by using linear combination of the original smoothness indicators to calculate the weights, the new scheme (WENO-Z) obtains superior results at almost the same computational cost of the WENO-JS scheme.

The WENO-Z method provides a straightforward way for improving the performance of a WENO scheme. Castro et al. \cite{Castro2011wenoz} developed higher-order WENO-Z schemes. Ha et al. \cite{Ha2013} constructed a new formula for local smoothness indicator and devised a new sixth-order global smoothness indicator. Fan et al. \cite{Faneta8} constructed several global smoothness indicators with truncation errors of up to eighth-order, the corresponding WENO-Z$\eta$ scheme can present fifth-order convergence in smooth regions, especially at critical points where the first and second derivatives vanish. Hu et al. \cite{Huweno6} constructed an adaptive central-upwind WENO-Z type scheme (WENO-CU6) in which a tunable parameter was introduced for the weighting function of WENO-Z to increase the contribution of optimal weights. Acker et al. \cite{Acker2016} presented a way of increasing the relevance of less smooth substencils by adding a new term into the WENO-Z weights; the new scheme (WENO-Z+) improves the resolution of high-frequency smooth waves. Liu et al. \cite{Liu2018} proposed a new sixth-order GSI (GSI-6) as well as a function consisted of the GSI-6 and local smoothness indicators (LSIs, $IS_k$) to calculate the weights of the fifth-order WENO scheme (WENO-ZA). For a smooth solution, the new weights satisfy the sufficient condition for fifth-order convergence in smooth regions; hence, the WENO-ZA scheme performs better than the WENO-Z scheme. For discontinuous solutions, the weights of WENO-ZA assigned to discontinuous substencils are as large as the ones of WENO-Z, i.e., the behavior of the WENO-ZA scheme in shock regions is similar to that of the WENO-Z scheme.

{For the nonlinear weight formulation $\alpha_k=c_k(1+(\frac{\tau_5}{IS_k+\epsilon})^q)$  of WENO-Z, Don and Borges \cite{Don2013} discussed the roles of two free parameters power  $\epsilon$, which is introduced to avoid zero denominator in $\alpha_k$, and $q$, which controls the amount of numerical dissipation. They proved that the optimal order of the WENO-Z scheme can be guaranteed with a much weaker condition $\epsilon=\Omega(\Delta x^m)$, where $m(r,q)\ge 2$ is the optimal sensitivity order, regardless of critical points. Recently, a modified nonlinear term $\Gamma=\Phi(\frac{\tau_5}{IS_k+\epsilon})^q$  was suggested by Wang et al.\cite{Wang_don_2019} to replace the original one $(\frac{\tau_5}{IS_k+\epsilon})^q$, where $\Phi$ is a function of a linear combination of the smoothness indicators $IS_k$. And an optimal variable $\epsilon=\Delta x^4$ with $q = 2$ was also suggested. Two new schemes are referred as WENO-D and WENO-A, they satisfy the Cp-property\cite{Wang_don_2019}.}

{For designing a WENO scheme, the ENO property should be put on the first place, and then, all operations should be in agreement with the physical requirements. Hence, although WENO-D/A can work well for those examples in Ref. \cite{Wang_don_2019}, its constructing method has some potential problems. First, the dimensions of length ($\Delta x$) and function $f$ (or $IS_k$) are different, so one can not simply put them together ($IS_k+\Delta x^m$). Second, if the modifier function $\Phi$, which has certain dimension related to the formula of $\Phi=min(1,\phi)$(here, one can not compare two variables with different dimensions, since $\phi$ has the dimension of $f$), is introduced, then one also can not do the addition operation of the linear term ($1$) and the nonlinear term, i.e., ($1+\Phi[\tau_5/(IS_k+\epsilon)]^q)$. Hence, the method may result in two issues, (1) the numerical solutions of WENO-D/A lose self-similarity, if different reference values are chosen to nondimensionalize the function $f$ and the computational region; (2) the resulted WENO-D/A schemes may lose the ENO property.}

{In order to obtain more accurate solution at critical points and avoid unmatched dimensions, in this paper, we propose a new method to construct a robust and high accurate WENO-Z type scheme. First, the square of the approximation of the fourth-order derivative, which is the maximal-order derivative can be approximated on a five-point stencil (the global stencil) by a fifth-order WENO scheme, is taken as the global smoothness indicator. Then, the constant $1$ is replaced by a function of the local smoothness indicators of the candidate sub-stencils. The function adaptively approaches to a small value if the global stencil contains a discontinuity and approaches to a large value if the global stencil is sufficiently smooth.}

This article is organized as follows: Section 2 describes the reconstruction procedure of several kinds of WENO schemes. Section 3 presents the new method for constructing high performance fifth-order WENO-Z type (WENO-ZN) scheme. Numerical experiments including one- and two-dimensional benchmark problems are presented in Section 4. Conclusions are drawn in Section 5.

\section{The fifth-order WENO schemes}

The one-dimensional scalar conservative law equation is used as a model to describe a numerical method
\begin{equation}\label{eq:1}
\frac{\partial{u}}{\partial{t}}+\frac{\partial{f(u)}}{\partial{x}}=0.
\end{equation}
The flux function $f(u)$ can be split into two parts as $f(u) = f^+(u) + f^-(u)$ with $df^+(u)/du\ge 0$ and $df^-(u)/du\le 0$. By defining the points $x_i=i\Delta{x}$, ($i=0,\dots,N$), where $\Delta{x}$ is the uniform grid spacing, the semi-discrete form of Eq.\eqref{eq:1} can be written as
\begin{equation}\label{eq:2}
\frac{du_i}{dt}=-\frac{\hat{f}_{i+1/2}-\hat{f}_{i-1/2}}{\Delta{x}},
\end{equation}
where $\hat{f}_{i\pm1/2}=\hat{f}^+_{i\pm1/2}+\hat{f}^-_{i\pm1/2}$ is the numerical flux. In this paper, only the positive part $\hat{f}^+_{i+1/2}$ is described and the superscript
$'+'$ is dropped for simplicity. The flux $\hat{f}^-_{i+1/2}$ is evaluated following the symmetric rule about $x_{i+1/2}$.

\subsection{The WENO-JS scheme}
The flux of a fifth-order WENO \cite{Jiang1996} scheme can be written as
\begin{equation}\label{eq:3}
\hat{f}_{i+1/2}=\sum_{k=0}^{2} \omega_kq_k,
\end{equation}
where $q_k$ is the third-order flux on the sub-stencil $S_{k}^3=(i+k-2,i+k-1,i+k)$, and given by
\begin{equation}\label{eq:4}
\begin{cases}
q_0=&\dfrac{1}{3}f_{i-2}-\dfrac{7}{6}f_{i-1}+\dfrac{11}{6}f_i, \\
q_1=&-\dfrac{1}{6}f_{i-1}+\dfrac{5}{6}f_i+\dfrac{1}{3}f_{i+1}, \\
q_2=&\dfrac{1}{3}f_i+\dfrac{5}{6}f_{i+1}-\dfrac{1}{6}f_{i+2} .
\end{cases}
\end{equation}
The weights $\omega_k$ of Jiang and Shu \cite{Jiang1996} is calculated as
\begin{equation}\label{eq:5}
\omega_k=\frac{\alpha_k}{\alpha_0+\alpha_1+\alpha_2},\\
\alpha_k=\frac{c_k}{(IS_k+\epsilon)^2},k=0,1,2,
\end{equation}
where, $IS_k$ is called as the local smoothness indicator (LSI), which is used to measure the relative smoothness of a solution on the sub-stencil $S_k$. $c_0=0.1,c_1=0.6$ and $c_2=0.3$ are the optimal weights, which generate the fifth-order upstream scheme. The parameter $\epsilon$ is a positive real number introduced to avoid the denominator becoming zero, and $\epsilon=10^{-6}$ is suggested by Jiang and Shu\cite{Jiang1996}.

In \cite{Jiang1996}, Jiang and Shu proposed a classical local smoothness indicator (LSI) as
\begin{equation}\label{eq:jiang}
IS_k=\sum_{l=1}^{r-1} \int_{x_{i-1/2}}^{x_{i+1/2}} (\Delta{x})^{2l-1}(q_k^{(l)})^2 dx,
\end{equation}
where, $q_k^{(l)}$ is the $l$th order derivative of $q_k(x)$, and $q_k(x)$ is the interpolation polynomial on sub-stencil $S_k^3$.

Taylor expansion of \eqref{eq:jiang} gives
\begin{equation}\label{eq:isjiang}
\begin{cases}
IS_0=\dfrac{13}{12}(f_{i-2}-2f_{i-1}+f_i)^2+\dfrac{1}{4}(f_{i-2}-4f_{i-1}+3f_i)^2 \\
IS_1=\dfrac{13}{12}(f_{i-1}-2f_{i}+f_{i+1})^2+\dfrac{1}{4}(f_{i-1}-f_{i+1})^2 \\
IS_2=\dfrac{13}{12}(f_{i}-2f_{i+1}+f_{i+2})^2+\dfrac{1}{4}(3f_{i}-4f_{i+1}+f_{i+2})^2
\end{cases}
\end{equation}

The Taylor expansion of $IS_k$ at $x_i$ for a smooth solution is often used to analyze the performance of a WENO scheme,
\begin{equation}\label{eq:7}
\begin{cases}
IS_0=&f_{i}'^2\Delta{x}^2+(\dfrac{13}{12}f_{i}''^2-\dfrac{2}{3}f_{i}'f_{i}''')\Delta{x}^4+
      (-\dfrac{13}{6}f_{i}''f_{i}'''+
      \dfrac{1}{2}f_{i}'f_{i}^{(4)})\Delta{x}^5+O(\Delta{x}^6), \\
IS_1=&f_{i}'^2\Delta{x}^2+(\dfrac{13}{12}f_{i}''^2+\dfrac{1}{3}f_{i}'f_{i}''')\Delta{x}^4+O(\Delta{x}^6), \\
IS_2=&f_{i}'^2\Delta{x}^2+(\dfrac{13}{12}f_{i}''^2-\dfrac{2}{3}f_{i}'f_{i}''')\Delta{x}^4+
      (\dfrac{13}{6}f_{i}''f_{i}'''-
      \dfrac{1}{2}f_{i}'f_{i}^{(4)})\Delta{x}^5+O(\Delta{x}^6) .
\end{cases}
\end{equation}

In \cite{Henrick2005}, Henrick et al. derived the necessary and sufficient conditions for fifth-order convergence of a fifth-order WENO scheme,
\begin{equation}\label{eq:ne_su}
\begin{cases}
\displaystyle\sum_{k=0}^{2}A_k(\omega^+_k-\omega^-_k)=O(\Delta x^3)\\
\omega^\pm-c_k=O(\Delta x^2)
\end{cases}
\end{equation}
where, $A_k$ are the coefficients of those terms with $\Delta x^3$ of the Taylor series expansions of $q_k$ (Eq.\eqref{eq:4}), $\omega^\pm$ are the weights of $\hat{f}_{i\pm 1/2}$ respectively.

Henrick et al. pointed out that the WENO-JS scheme may even decrease to third-order accuracy at critical points, hence a mapping function \cite{Henrick2005} is proposed to make the new weights satisfy a sufficient condition, which is given as
\begin{equation}\label{eq:10}
\omega^\pm_k-c_k=O(\Delta{x}^3).
\end{equation}
Although this condition is not necessary, as mentioned by Henrick et al. \cite{Henrick2005}, Eq.\eqref{eq:10} can serve as a simple criteria to design the weights for fifth-order WENO schemes.

\subsection{The WENO-Z scheme}

The fifth-order WENO-Z scheme is proposed by Borges et al. \cite{Borges2008} by introducing a global smoothness indicator (GSI) $\tau_5$ to calculate the weights,
\begin{equation}\label{eq:15}
\omega_k^{Z}=\frac{\alpha_k}{\alpha_0+\alpha_1+\alpha_2},\
\alpha_k=c_k\left(1+(\frac{\tau_5}{IS_k+\epsilon})^q\right) .
\end{equation}
The original $\tau_5$ of Borges et al. is
\begin{equation}\label{eq:tao5}
\tau_5=\left|IS_2-IS_0\right| \  .
\end{equation}
Applying the Taylor expansions of $IS_k$ \eqref{eq:7}, there is
\begin{equation}\label{eq:15a}
\begin{split}
\tau=\left|\frac{13}{3}f''f'''-f'f^{(4)}\right|\Delta{x}^5+O(\Delta{x}^6).
\end{split}
\end{equation}
Hence, one can get
\begin{align}\label{eq:16}
\omega_k^{\tau_5}=
\begin{cases}
c_k+O(\Delta{x}^{3q}),& f_{i}'\ne0, \\
c_k+O(\Delta{x}^q),& f_{i}'=0,
\end{cases}&
\end{align}
where $q$ is a tunable parameter. Numerical results in \cite{Borges2008} demonstrated that, if $q$ takes $1$, the accuracy order at critical points is only fourth; with $q=2$, the scheme can achieve fifth-order accuracy. Meanwhile, Borges et al. pointed out that, for solutions containing discontinuities, increasing $q$ makes the scheme more dissipative. As lower dissipation of WENO-Z is more desirable than its rate of convergence at critical points when solving problems involving shocks, $q=1$ is suggested for the WENO-Z scheme in \cite{Borges2008}.

\subsection{Several improved WENO-Z-type schemes}\label{num_meth}
The weight function \eqref{eq:15} of the WENO-Z scheme provides a straight-forward guideline for improving the accuracy of a WENO scheme. For completeness, here, several improved WENO-Z-type schemes are briefly introduced (please refer to \cite{Liu2018} for more details).

(1) WENO-NS: Ha et al.\cite{Ha2013} constructed a sixth-order global smoothness indicator (GSI) as
\begin{equation} \label{eq:ha}
\zeta=\frac{1}{2}\left(|\beta_0-\beta_2|^2+g(|L_{1,1}f|)^2\right),
\end{equation}
where, $\beta_k=\xi|L_{1,k}f|+|L_{2,k}f|$,  $\xi$ is a tunable parameter that governs the tradeoff between the accuracies around smooth region and discontinuous region, and $L_{l,k}f$ is the approximation of the $l$th derivative $f^{(l)}_{i+1/2}\Delta x^l$ on sub-stencil $S_k$.
The local smoothness indicator (LSI) $IS_k$ is calculated as $IS_k=\beta_k^2$.

(2) WENO-P: Kim et al. \cite{Kim2016} simplified the sixth-order GSI Eq.\eqref{eq:ha} as
\begin{equation} \label{eq:hap}
\zeta=(\beta_0-\beta_2)^2,
\end{equation}
to reduce computation cost, and introduced a parameter $\delta$ to make a balanced contribution of the $\beta_k$ of Ha et al. as
\begin{equation} \label{eq:kim}
\tilde{\beta_0}=\beta_0, \ \tilde{\beta_1}=(1+\delta)\beta_1,\ \tilde{\beta_2}=(1-\delta)\beta_2.
\end{equation}

(3) WENO-$\eta$: Fan et al. \cite{Faneta8} proposed a sixth-order and two eighth-order GSIs as
\begin{equation}\label{eq:eta8}
\begin{split}
&\tau_6=|\eta_5-\frac{IS_0+4IS_1+IS_2}{6}|,\\
&\tau_{81}=|(|f_0^{(1)}|-|f_2^{(1)}|)(f_0^{(2)}+f_2^{(2)}-2f_1^{(2)})|, \\
&\tau_{82}=(|f_0^{(1)}|-|f_2^{(1)}|)^2+(f_0^{(2)}+f_2^{(2)}-2f_1^{(2)})^2,
\end{split}
\end{equation}
where, $\eta_5=\dfrac{1}{144}{[(f_{i-2}-8f_{i-1}+8f_{i+1}-f_{i+2})^2+(f_{i-2}-16f_{i-1}+30f_{i}-16f_{i+1}+f_{i+2})^2]}$, and the local smoothness indicator $IS_k$ takes the formula suggested by Shen and Zha \cite{ShenAIAA2008},
\begin{equation}\label{eq:iszha}
\begin{split}
IS_k=\sum_{l=1}^{r-1}\gamma_l\Delta x^{2l}[f_k^{(l)}]^2
\end{split}
\end{equation}
where $f_k^{(l)}$ is the approximation of the $l$th order derivative $f^{(l)}_{x_i}$ on sub-stencil $S_k$. The application of formula \eqref{eq:iszha} is flexible and convenient, for example, Jiang and Shu's formula\eqref{eq:isjiang} gives $\gamma_1=1$ and $\gamma_2=13/12$ while Fan et al. chose $\gamma_1=1$ and $\gamma_2=1$.

(4) WENO-CU6: Hu et al.\cite{Huweno6} constructed an adaptive central-upwind WENO scheme, in which, the weights are constructed as
\begin{equation}\label{eq:hu}
\omega_k=\frac{\alpha_k}{\alpha_0+\cdots +\alpha_3},\
\alpha_k=c_k\left(C+\frac{\tau_6}{IS_k+\epsilon}\right) ,
\end{equation}
where, $IS_3$ is the smoothness indicator \eqref{eq:jiang} on the stencil $S^6=(i-2,\cdots,i+3)$, and $\tau_6=|IS_3-(IS_0+4IS_1+IS_2)/6|$. The parameter $C$  is introduced to increase the contribution of optimal weights and decrease numerical dissipation, and $C=20$ is suggested in \cite{Huweno6}.

(5) WENO-Z+: Acker et al. \cite{Acker2016} proposed a way of improving the results of WENO-Z by increasing the weights of less-smooth sub-stencils,
\begin{equation}\label{eq:ac}
\alpha_k=c_k\left[1+\left(\frac{\tau_5+\epsilon}{IS_k+\epsilon}\right)^2+\lambda\left(\frac{IS_k+\epsilon}{\tau_5+\epsilon}\right)\right],
\end{equation}
where $\lambda$ is a parameter being dependent on the grid spacing.


(6) WENO-ZA: Liu et al. \cite{Liu2018} proposed a new method to calculate the weights,
\begin{equation}\label{eq:nw}
\omega_k=\frac{\alpha_k}{\alpha_0+\alpha_1+\alpha_2},\
\alpha_k=c_k\left(1+A\frac{\tau}{IS_k+\epsilon}\right),
\end{equation}
where, the function $A$ is
\begin{equation}\label{eq:np}
A=\frac{\tau}{IS_0+IS_2-\tau+\epsilon}.
\end{equation}
and the global smoothness indicator (GSI) is designed as
\begin{equation}\label{eq:ntao}
\begin{split}
\tau=\gamma_1(|f_0^{(1)}|-|f_2^{(1)}|)^2+\gamma_2(|f_0^{(2)}|-|f_2^{(2)}|)^2.
\end{split}
\end{equation}

{(7) WENO-D/A: Recently, Wang et al.\cite{Wang_don_2019} constructed the WENO-D/A schemes. The un-normalized weight $\alpha_k$ of WENO-D is
\begin{equation}\label{eq:weno_d}
\alpha_k=c_k\left(1+\Phi(\frac{\tau_5}{IS_k+\epsilon})^q\right),\
\end{equation}
where,
$$\Phi=min(1,\phi), \text{and}\ \phi=\sqrt{|IS_0-2IS_1+IS_2}$$
WENO-D has the similar form as WENO-ZA\cite{Liu2018}.}

{WENO-A is a modification of WENO-D, its weight is
\begin{equation}\label{eq:weno_a}
\alpha_k=c_k(max(1,\Phi(\frac{\tau_5}{IS_k+\epsilon})^q)\
\end{equation}}

{Wang et al.\cite{Wang_don_2019} analyzed that the WENO-D/A schemes satisfy the Cp-property. However, since $\phi$ has the dimension of $f$, one cannot simply compare constant $1$ and $\phi$, such as minimal function and maximum function, and also cannot do addition operation of $(1+\Phi(\frac{\tau_5}{IS_k+\epsilon})^q)$. Otherwise, if different reference values are chosen to nondimensionalize the function $f$, the numerical solutions of WENO-D/A lose self-similarity. In addition, if a large reference value is used, the resulted WENO-D/A schemes may generate oscillation. Since flux $f$ and length (or $\Delta x$) have different dimensions, the similar issues also exist if $\epsilon$ takes a function of the grid spacing $\Delta x$, such as $\epsilon=\Delta x^4$ suggested in \cite{Wang_don_2019}. Spurious numerical solutions caused by these issues will be numerically demonstrated in Sec.\ref{num_tests}.}

\section{The new WENO scheme}

In this section, we propose a new method to calculate the weights $\alpha_k$ for a WENO-Z type scheme, based on the following analysis. First, from the formulation of $\alpha_k$ (Eq.\eqref{eq:15}), there are four parameters independent of $k$, i.e., the constant $1$, $q$, $\epsilon$, and $\tau$. There are many papers \cite{Jiang1996, Henrick2005, Borges2008,Don2013, ShenAIAA2008} discussing the roles of the two parameters $\epsilon$ and $q$. Here, we take into account the constant $1$ and GSI $\tau$. It is clear that, the constant $1$ can be replaced by a function independent of $k$, the function is required to approach to a large value for a smooth global stencil $S^5$ for obtaining low dissipation and high accuracy; at the same time, if the global stencil $S^5$ contains a discontinuity, the function is required to approach to a value small enough to keep the ENO property. Then, since derivatives of any order can reflect the information of discontinuous solution to a certain extent, we can use the derivative of the highest order that can be approximated on the global stencil as the global smoothness indicator to achieve the maximal order of accuracy for smooth solution. For this purpose, we suggest a function as
\begin{equation}\label{eq:function}
C=A\left(\frac{IS_0+IS_2-\tau+\epsilon}{\tau+\epsilon}\right)^2
\end{equation}
where, $A$ is a constant, in this paper we take $A=10$ and this choice will be discussed later, $\tau$ can take $\tau_5$ suggested by Borges et al.\cite{Borges2008}, i.e.,
\begin{equation}\label{eq:tau}
\tau=\tau_5=|IS_0-IS_2|
\end{equation}
and
\begin{equation}\label{eq:tau8}
\tau_8=(f_{i-2}-4f_{i-1}+6f_i-4f_{i+1}+f_{i+2})^2
\end{equation}
is used as the global smoothness indicator.

The un-normalized weight $\alpha_k$ is then calculated by
\begin{equation}\label{eq:new}
\alpha_k=c_k\left(C+\frac{\tau_8}{IS_k+\epsilon}\right) .
\end{equation}
For convenience, we call the resulted scheme with the new weight \eqref{eq:new} as the WENO-ZN scheme.

Now, we discuss the properties of the new weight:

(1) For a smooth global stencil, the Taylor series expansion of \eqref{eq:tau8} gives
\begin{equation}
\tau_8=\left(f^{(4)}_i\Delta x^4+O(\Delta x^5)\right)^2
\end{equation}
By using the formula \eqref{eq:new}, there is
\begin{align}\label{eq:C1}
\alpha_k=
\begin{cases}
c_k\left(C+O(\Delta x^4)\right),& \text{if}\ f_{i}'=0, \\
c_k\left(C+O(\Delta x^2)\right),& \text{if}\ f_{i}'=0\ \text{and}\ f_{i}''=0
\end{cases}&
\end{align}
Hence, the new weights can satisfy the sufficient condition for fifth-order convergence \eqref{eq:10} at first critical point ($f'_i=0$), and can achieve fourth-order accuracy at second-order critical point ($f'_i=0$ and $f''_i=0$).

Meanwhile, from the Taylor series expansion \eqref{eq:7}, $IS_0$ and $IS_2$ always have the same first term, hence there is
\begin{align}\label{eq:C}
C=
\begin{cases}
O(\Delta{x}^{-6}),& \text{if}\ f_{i}'\ne0, \\
O(\Delta{x}^{-2}),& \text{otherwise}
\end{cases}&
\end{align}
Eq.(\ref{eq:C}) means that $C$ is a large value($>>1$).

As the value of $\tau_8$ is independent of $k$, by setting $C$ to a large value, the influence of the difference between $IS_k$s on the un-normalized weights can be overwhelmed and the numerical dissipation of the resulted scheme is therefore reduced.

(2) If the solution is discontinuous on a global stencil, without loss of generality, assume that the sub-stencil $S_0^3$ contains a discontinuity and $S_2^3$ is smooth, then there is $IS_0 \gg IS_2$, and
\begin{equation} \label{eq:nouse1}
\tau \approx IS_0 \gg IS_2\
\end{equation}
and hence $C\ll 1$, the relative magnitudes of $\alpha_k(k=0,1,2)$ are mainly determined by the second term in \eqref{eq:new}. Since $\tau_8$ is independent of $k$, the ENO property of the resulted scheme can be kept well by the local smoothness indicator $IS_k$.

\begin{figure}
\begin{minipage}{0.5\linewidth}
\centering
\includegraphics[width=1.0\textwidth]{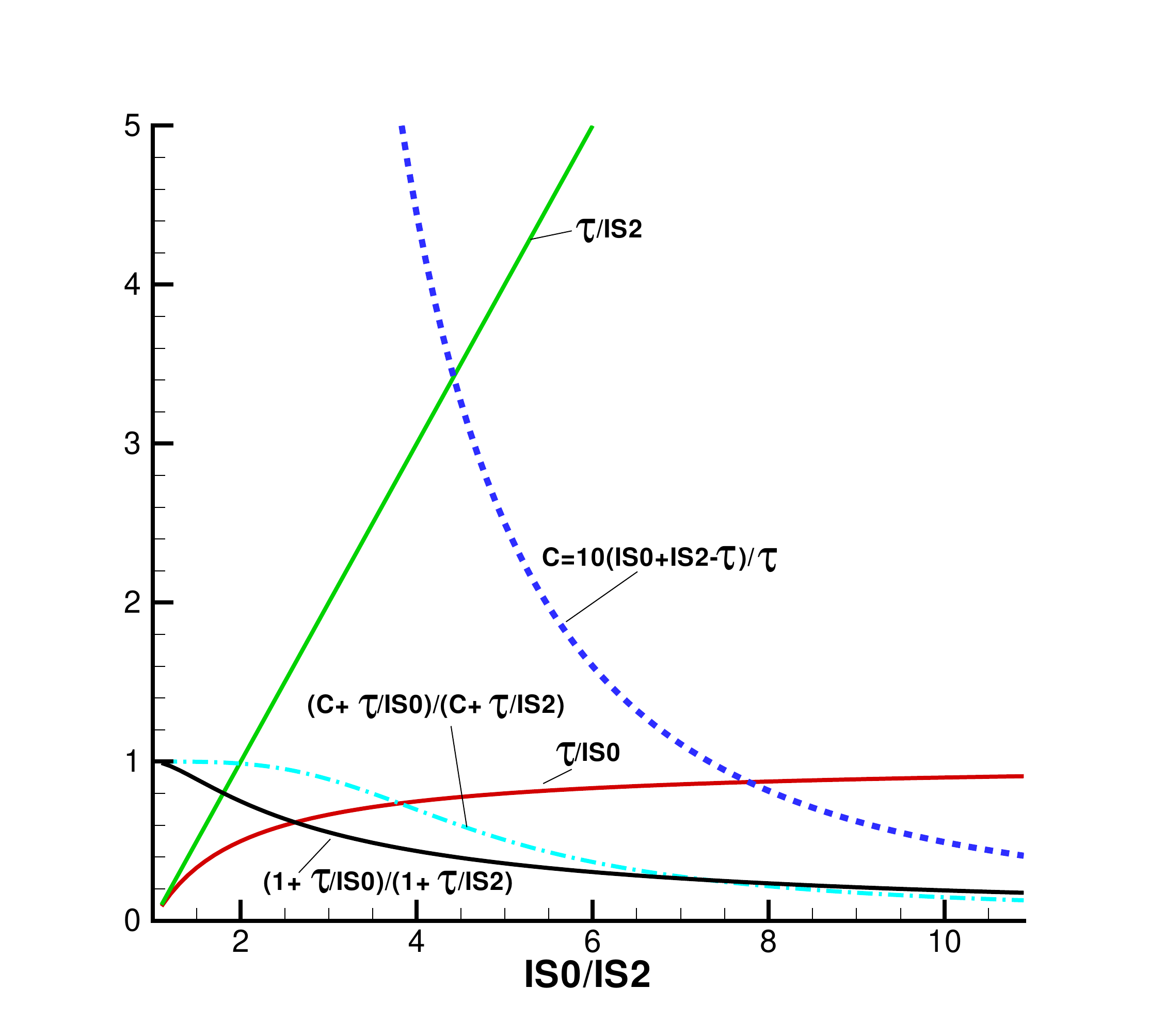}
\caption{Various distribution vs $IS0/IS2$ }
\label{fig:c-is0-is2_01}
\end{minipage}
\begin{minipage}{0.5\linewidth}
\centering
\includegraphics[width=1.0\textwidth]{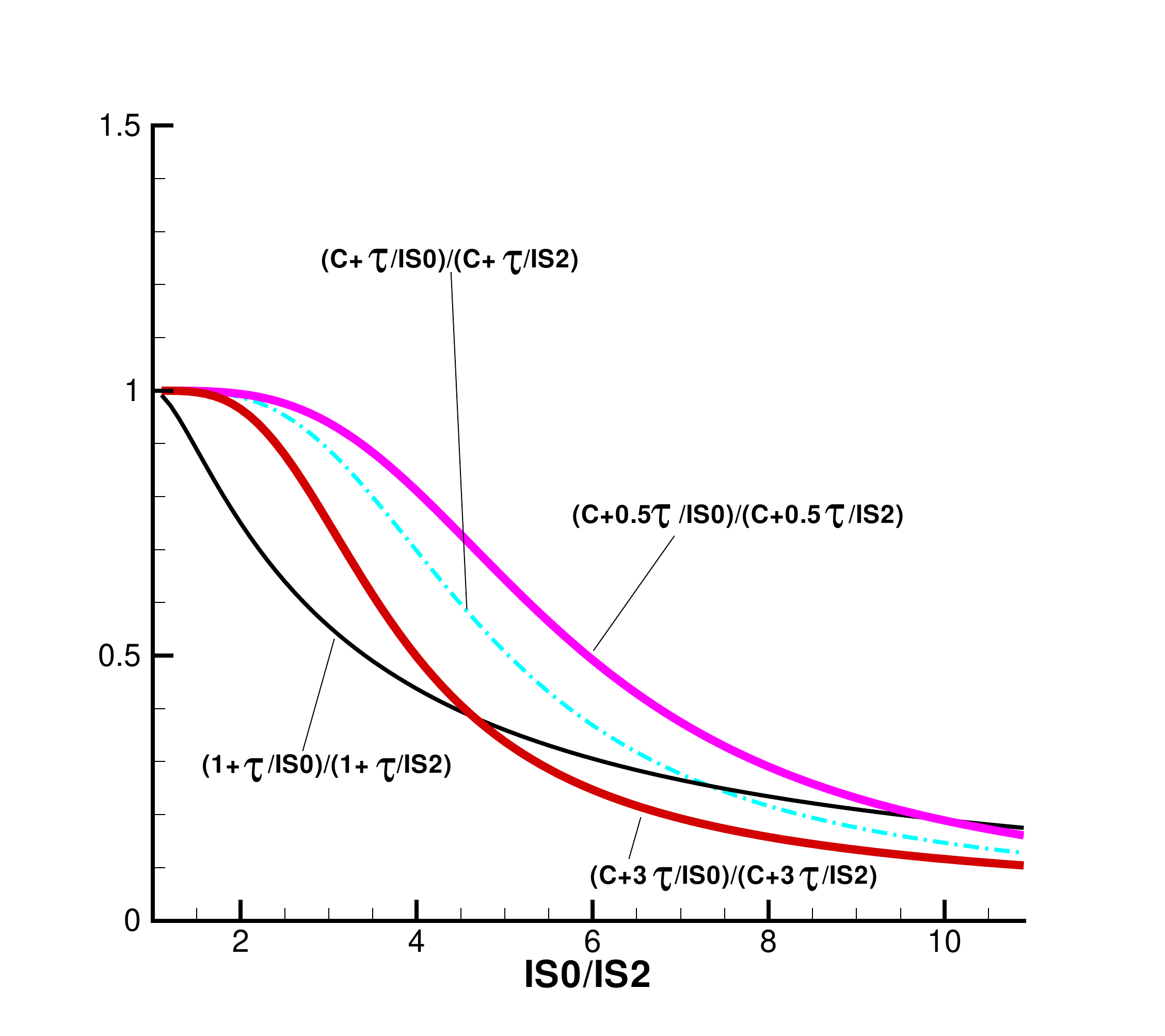}
\caption{Various distribution vs $IS0/IS2$}
\label{fig:c-is0-is2_02}
\end{minipage}
\end{figure}

(3) Without loss of generality, the following formula can be used to discuss the original method of WENO-Z and the new weighting method,
  $$(1+\frac{\tau}{IS_k+\epsilon})\ \ \ and \ \  (C+\frac{\tau}{IS_k+\epsilon}),$$
and the influence of $\epsilon$ is neglected. Fig.\ref{fig:c-is0-is2_01} shows various distributions vs $(IS_0/IS_2)$.
From this figure, it can be seen that, in a large range of $(IS_0/IS_2)$, the new method gives a more balanced contribution (the ratio $[C+\tau/S_0]/[C+\tau/IS_2]\sim 1$) of $IS_0$ and $IS_2$ than the original method does, hence the new method is less dissipative. If $(IS_0/IS_2)$ is large enough, the contribution of $IS_0$ ($S_0^3$ is regarded as the discontinuous sub-stencil) in the new method is less than that in the original one, this is helpful to increase the shock-capturing capability (ENO property). This property is also used to decide the constant $A$ in \eqref{eq:function}, i.e., if $(IS_0/IS_2)$ is larger than one order of magnitude ($IS_0/IS_2>10$), then the new method can satisfy the requirement that the contribution of $IS_0$ is not larger than that in the original WENO-Z scheme.

Table \ref{tau_58} gives the coefficient of $[f_{i+k}-f_{i+k-1}]^2(k=-1,\cdots, 2)$ in $\tau_5$ and $\tau_8$. Theoretically, if there only exists one discontinuity at the global stencil $S^5$, Table \ref{tau_58} indicates that $\tau_8$ is almost $0.75\sim 2.7$ times of $\tau_5$. In Fig.\ref{fig:c-is0-is2_02}, the two curves of $[C+0.5\tau/S_0]/[C+0.5\tau/IS_2]$) and $[C+3\tau/S_0]/[C+3\tau/IS_2]$) are also plotted. Clearly, the two curves both meet the requirement above, hence $A=10$ is reasonable for the fifth-order scheme.

\begin{table}
  \caption{The coefficient of $(f_{i+k}-f_{i+k-1})^2, k=-1,\cdots, 2$}\label{tau_58}
  \begin{center}\footnotesize
  \begin{tabular}{cccccc}
  \toprule
    & $(f_{i-1}-f_{i-2})^2$ & $(f_{i}-f_{i-1})^2$ & $(f_{i+1}-f_{i})^2$ &  $(f_{i+2}-f_{i+1})^2$ \\
  \midrule
  $\tau_5$ &  4/3 &  10/3 &  10/3 &  4/3 \\
  $\tau_8$  &  1 &  9 &  9 &  1 \\
   \bottomrule
  \end{tabular}
  \end{center}
\end{table}

(4) The following function $u_0(x)$ with a discontinuous point $x=0$ is used to show that $\tau_8$ \eqref{eq:tau8} has a similar behaviour as $\tau_5$ \eqref{eq:tao5}, and hence $\tau_8$ can be used as a global smoothness indicator.
\begin{align}\label{eq:single}
u_0(x)=&
\begin{cases}
-\text{sin}(\pi x)-\dfrac{1}{2}x^3,&-1\le x< 0, \\
-\text{sin}(\pi x)-\dfrac{1}{2}x^3+1,&0\le x \le 1.
\end{cases}&
\end{align}

Table \ref{function_1} gives the numerical results of the function $u_0(x)$. In this table, for the cases of $IS_0\ge IS_2$,
  $$
  R=\frac{1+\tau_5/(IS_0+\epsilon)}{1+\tau_5/(IS_2+\epsilon)} \ \ \ and \ \ \ R'=\frac{C+\tau_8/(IS_0+\epsilon)}{C+\tau_8/(IS_2+\epsilon)},
  $$
are used to measure the contributions of the sub-stencils $S_0^{3}$ and $S_2^{3}$ , while for the cases of $IS_2>IS_0$, those values are calculated by
$$
R=\frac{1+\tau_5/(IS_2+\epsilon)}{1+\tau_5/(IS_0+\epsilon)} \ \ \ and \ \ \ R'=\frac{C+\tau_8/(IS_2+\epsilon)}{C+\tau_8/(IS_0+\epsilon)},
$$
and a negative sign $'-'$ is assigned to them. Table \ref{function_1} shows that, the values of $\tau_5$ and $\tau_8$ at discontinuous (global) stencils are much larger than those at smooth stencils. In addition, at those discontinuous stencils, $\tau_8$ is almost $0.7\sim 3$ times of $\tau_5$, this is in agreement with the analysis above. It also can be seen that, the contribution of discontinuous sub-stencils in the new method is about half of that in original method of WENO-Z, i.e. $R'<R$ (neglecting the negative sign $'-'$).


\begin{table}
  \caption{The numerical results of the function $u_0(x)$ \eqref{eq:single}}\label{function_1}
  \begin{center}\footnotesize
  \begin{tabular}{cccccc}
  \toprule
   $x_i$ & $u_0(x_i)$ & $\tau_5$ & $\tau_8$ &  $R$ & $R'$ \\
  \midrule
  -0.6000E-01 &  0.1875E+00 &  0.4526E-06 &  0.8518E-11 &  0.1000E+01 &  0.1000E+01\\
  -0.4000E-01 &  0.1254E+00 &  0.3065E-06 &  0.3811E-11 &  0.1000E+01 &  0.1000E+01\\
  -0.2000E-01 &  0.6279E-01 &  0.1396E+01 &  0.1000E+01 & -0.5625E-02 & -0.2818E-02\\
   0.0000E-00 &  0.0000E-00 &  0.3144E+01 &  0.9000E+01 & -0.2513E-02 & -0.1257E-02\\
   0.2000E-01 &  0.9372E+00 &  0.3145E+01 &  0.9000E+01 &  0.2503E-02 &  0.1252E-02\\
   0.4000E-01 &  0.8746E+00 &  0.1395E+01 &  0.1000E+01 &  0.5569E-02 &  0.2790E-02\\
   0.6000E-01 &  0.8125E+00 &  0.4526E-06 &  0.8518E-11 &  0.1000E+01 &  0.1000E+01\\
   0.8000E-01 &  0.7511E+00 &  0.5903E-06 &  0.1500E-10 &  0.1000E+01 &  0.1000E+01\\
   \bottomrule
  \end{tabular}
  \end{center}
\end{table}

\section{Numerical examples}\label{num_tests}
In this section, several problems, including linear advection problems and one- and two-dimensional Euler problems, are considered to evaluate the performance of the new scheme. The time derivative is approximated with the third-order TVD Runge-Kutta method \cite{Shu1988}. Unless noted otherwise, the CFL number always takes 0.5 in this paper.

As we pointed out in Sec. \ref{num_meth}, the numerical results of WENO-D/A lose the self-similarity if different reference values are used to nondimensionalize the unknown variable (or computational region), and the numerical results may be oscillatory, so only several examples are calculated by WENO-D and used to demonstrate the drawbacks of WENO-D. In addition, various numerical results have shown that those improved WENO-Z-type schemes mentioned in Sec.2.3 perform well in most of the tested cases. However, those parameters, such as $\xi$ in \eqref{eq:ha} and \eqref{eq:hap}, $\delta$ in \eqref{eq:kim}, $C$ in \eqref{eq:hu}, and $\lambda$ in \eqref{eq:ac}, are problem-dependent and chosen empirically, and the WENO-Z$\eta$ schemes which use high order GSIs \eqref{eq:eta8}, are prone to generate oscillations near discontinuities. Since the numerical comparisons \cite{Liu2018} of WENO-ZA, WENO-Z, WENO-Z$\eta$, WENO-NS, and WENO-P showed comprehensive advantages (including ENO property, high-order accuracy, high resolution, and low dissipation) of WENO-ZA over the others, this paper only considers the comparisons of WENO-Z, WENO-ZA, and the present scheme.

\subsection{The accuracy at critical point}

{\renewcommand\baselinestretch{1.5}
\begin{table}
  \caption{Convergence order at the critical point}\label{convergence}
  \begin{center}\footnotesize
  \begin{tabular}{cccccccc}
  \toprule
  \multirow{2}{*}{Case} &
  \multirow{2}{*}{$\Delta{x}$} &
  \multicolumn{2}{c}{WENO-Z} &
  \multicolumn{2}{c}{WENO-ZA} &
  \multicolumn{2}{c}{Present} \\
  & & error & order & error & order & error & order \\
  \midrule
    &   0.2500E-01&    0.963652E-09&      ---&    0.964557E-09&      ---&    0.964557E-09&      ---\\
    &   0.1250E-01&    0.303249E-10&    4.990&    0.303284E-10&    4.991&    0.303284E-10&    4.991\\
    &   0.6250E-02&    0.950693E-12&    4.995&    0.950706E-12&    4.996&    0.950706E-12&    4.996\\
  $k=1$  &   0.3125E-02&    0.297558E-13&    4.998&    0.297559E-13&    4.998&    0.297559E-13&    4.998\\
    &   0.1563E-02&    0.930596E-15&    4.999&    0.930596E-15&    4.999&    0.930596E-15&    4.999\\
    &   0.7813E-03&    0.290925E-16&    4.999&    0.290925E-16&    4.999&    0.290925E-16&    4.999\\
    &   0.3906E-03&    0.909317E-18&    5.000&    0.909317E-18&    5.000&    0.909317E-18&    5.000\\
   &   0.1953E-03&    0.284189E-19&    5.000&    0.284189E-19&    5.000&    0.284189E-19&    5.000\\
   \midrule
    &   0.2500E-01&    0.124183E-05&      ---&    0.558942E-08&      ---&    0.481106E-08&      ---\\
    &   0.1250E-01&    0.724845E-07&    4.099&    0.156223E-09&    5.161&    0.151455E-09&    4.989\\
    &   0.6250E-02&    0.433711E-08&    4.063&    0.478277E-11&    5.030&    0.475058E-11&    4.995\\
  $k=2$  &   0.3125E-02&    0.264348E-09&    4.036&    0.148965E-12&    5.005&    0.148733E-12&    4.997\\
    &   0.1563E-02&    0.162992E-10&    4.020&    0.465399E-14&    5.000&    0.465225E-14&    4.999\\
    &   0.7813E-03&    0.101153E-11&    4.010&    0.145464E-15&    5.000&    0.145451E-15&    4.999\\
    &   0.3906E-03&    0.629935E-13&    4.005&    0.454651E-17&    5.000&    0.454641E-17&    5.000\\
   &   0.1953E-03&    0.392993E-14&    4.003&    0.142093E-18&    5.000&    0.142092E-18&    5.000\\
   \midrule
    &   0.2500E-01&    0.544997E-03&      ---&    0.288573E-03&      ---&    0.283929E-03&      ---\\
    &   0.1250E-01&    0.128406E-03&    2.086&    0.723469E-04&    1.996&    0.681519E-04&    2.059\\
    &   0.6250E-02&    0.293311E-04&    2.130&    0.181277E-04&    1.997&    0.145789E-04&    2.225\\
   $k=3$ &   0.3125E-02&    0.653022E-05&    2.167&    0.453707E-05&    1.998&    0.231768E-05&    2.653\\
    &   0.1563E-02&    0.144748E-05&    2.174&    0.113490E-05&    1.999&    0.237078E-06&    3.289\\
    &   0.7813E-03&    0.327355E-06&    2.145&    0.283803E-06&    2.000&    0.176821E-07&    3.745\\
    &   0.3906E-03&    0.763689E-07&    2.100&    0.709605E-07&    2.000&    0.116339E-08&    3.926\\
    &   0.1953E-03&    0.183126E-07&    2.060&    0.177413E-07&    2.000&    0.737522E-10&    3.980\\
   \bottomrule
  \end{tabular}
  \end{center}
\end{table}
}

The function $f(x)=x^kexp(x)$ is used to test the convergence rate of a WENO scheme at critical point \cite{Faneta8}. For the cases with $k=2$ and $k=3$, the point $x=0$ is a first-order critical point and a second-order critical point respectively.
As shown in Table \ref{convergence}, the original WENO-Z scheme with $q=1$ only gets fourth order accuracy and second order accuracy for the cases with first-order critical point and second-order critical point. The WENO-ZA scheme achieves fifth-order accuracy for the first-order critical point, but it is only second order for the second-order critical point. While, the present scheme even can reach fourth-order accuracy for the second-order critical point.

\subsection{Linear advection problems}
In the following, we test the accuracy of WENO schemes for the linear advection equation
\begin{align}\label{eq:26}
\begin{cases}
\dfrac{\partial{u}}{\partial{t}}+\dfrac{\partial{u}}{\partial{x}}=0,& x_0\leqslant x\leqslant{x_1},\\
u(x,t=0)=u_0(x),& \text{periodic boundary}.
\end{cases}
\end{align}
The exact solution of Eq.\eqref{eq:26} is given by
\begin{equation}
u(x,t)=u_0(x-t).
\end{equation}
Three linear cases are calculated.

Case 1: The initial condition is
\begin{align}\label{eq:1st_case}
u_0(x)=\frac{1}{A}\left\{
\begin{array}{ll}
\frac{1}{6}(G(x,\beta,z-\delta)+G(x,\beta,z+\delta)+4G(x,\beta,z)),& -0.8\le x\le -0.6,\\
1,& -0.4\le x \le -0.2,\\
1-|10(x-0.1)|,&\ 0\le x \le 0.2,\\
\frac{1}{6}(F(x,\alpha,\alpha-\delta)+F(x,\alpha,\alpha+\delta)+4F(x,\alpha,a)),&\ 0.4\le x\le 0.6,\\
0,& \ otherwise
\end{array}
\right .
\end{align}
where, $G(x,\beta,z)=e^{-\beta (x-z)^2}$, $F(x,\alpha,a)=\sqrt{max(1-\alpha^2(x-a)^2,0)}$. Same as in Ref.\cite{Jiang1996}, the constants for this case are taken as $a=0.5$, $z=-0.7$, $\delta=0.005$, $\alpha=10$, and $\beta=log2/36\delta^2$. The solution contains a smooth combination of Gaussians, a square wave, a sharp triangle wave, and a half ellipse. For this case, the computational time is $T=6$.

\begin{figure}
\begin{minipage}{0.5\linewidth}
\centering
\includegraphics[width=1.0\textwidth]{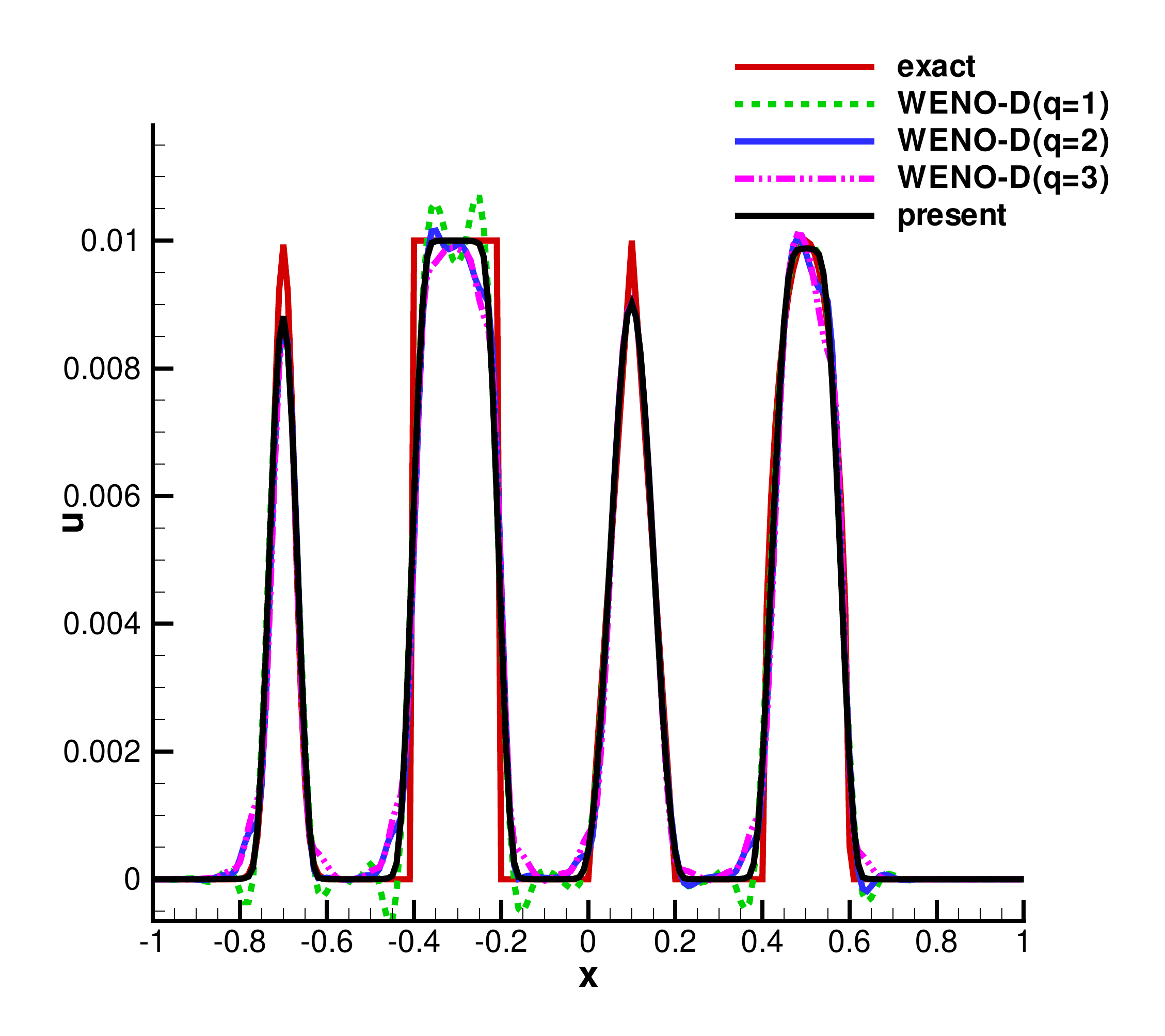}
\caption{Numerical results of case 1(\ref{eq:1st_case}) , A=100}
\label{fig:4wave_A100}
\end{minipage}
\begin{minipage}{0.5\linewidth}
\centering
\includegraphics[width=1.0\textwidth]{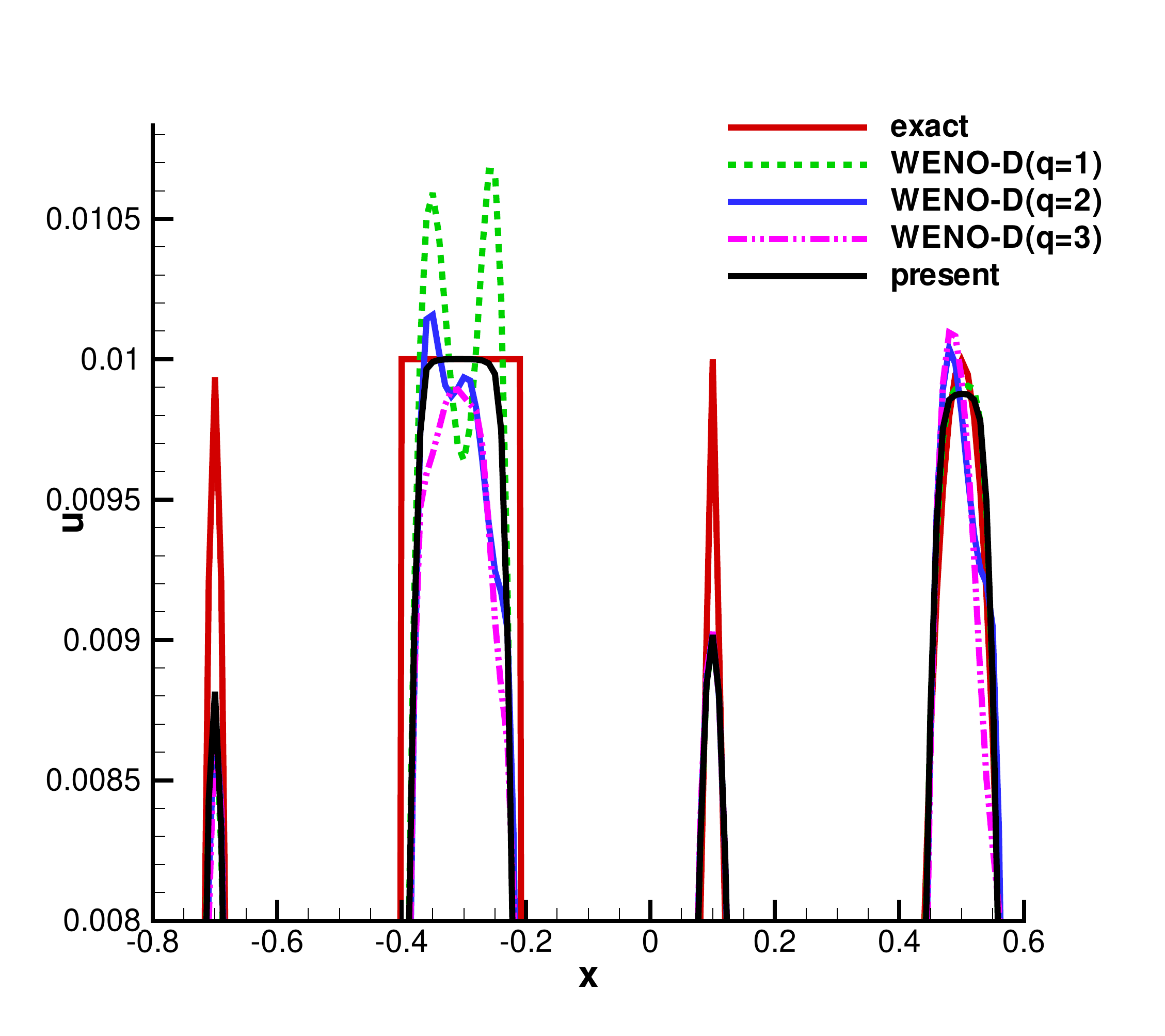}
\caption{Enlarged plot of Fig.\ref{fig:4wave_A100}}
\label{fig:4wave_A100_zoom}
\end{minipage}
\end{figure}

\begin{figure}
\begin{minipage}{0.5\linewidth}
\centering
\includegraphics[width=1.0\textwidth]{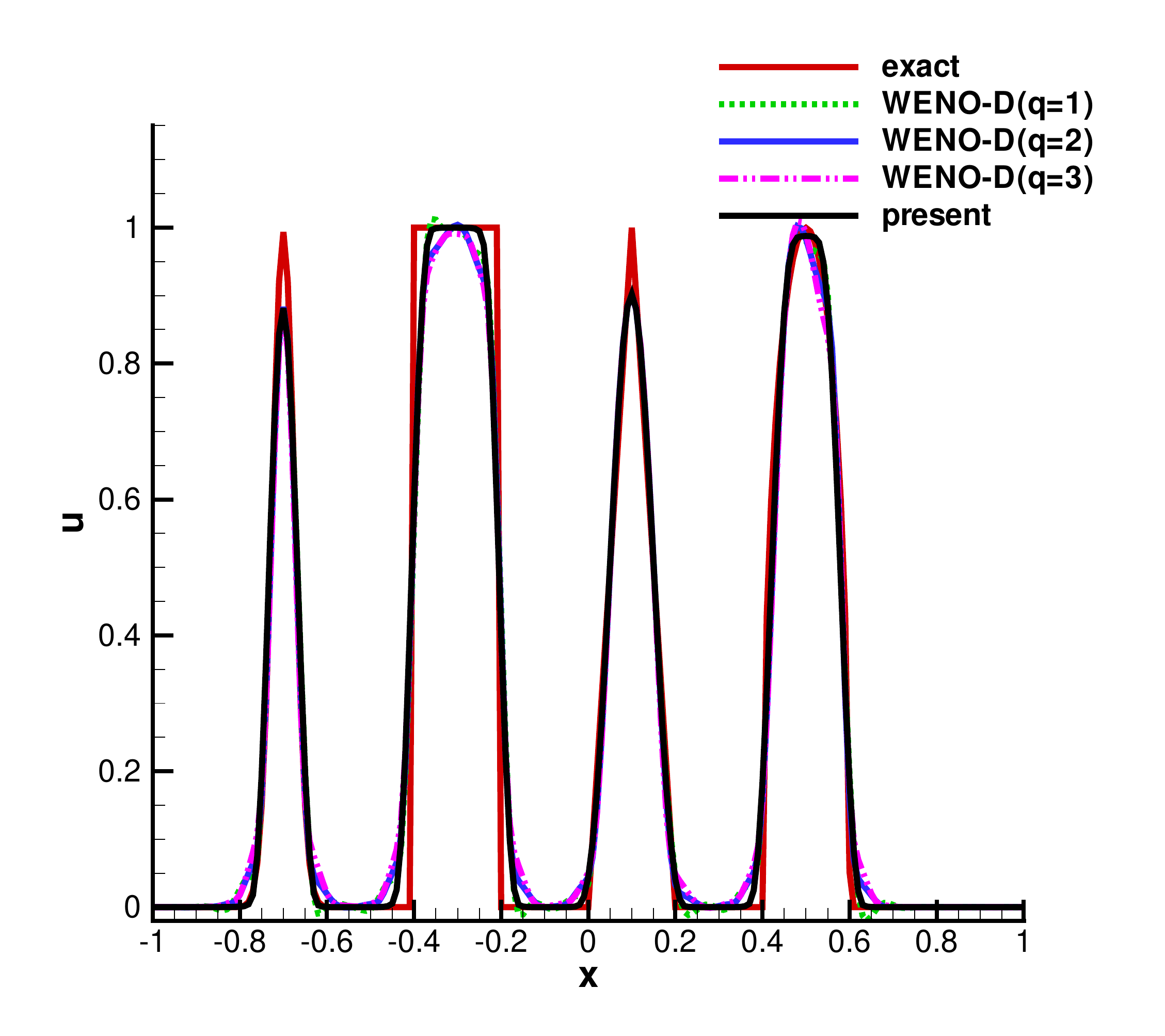}
\caption{Numerical results of case 1(\ref{eq:1st_case}), A=1}
\label{fig:4wave_A1}
\end{minipage}
\begin{minipage}{0.5\linewidth}
\centering
\includegraphics[width=1.0\textwidth]{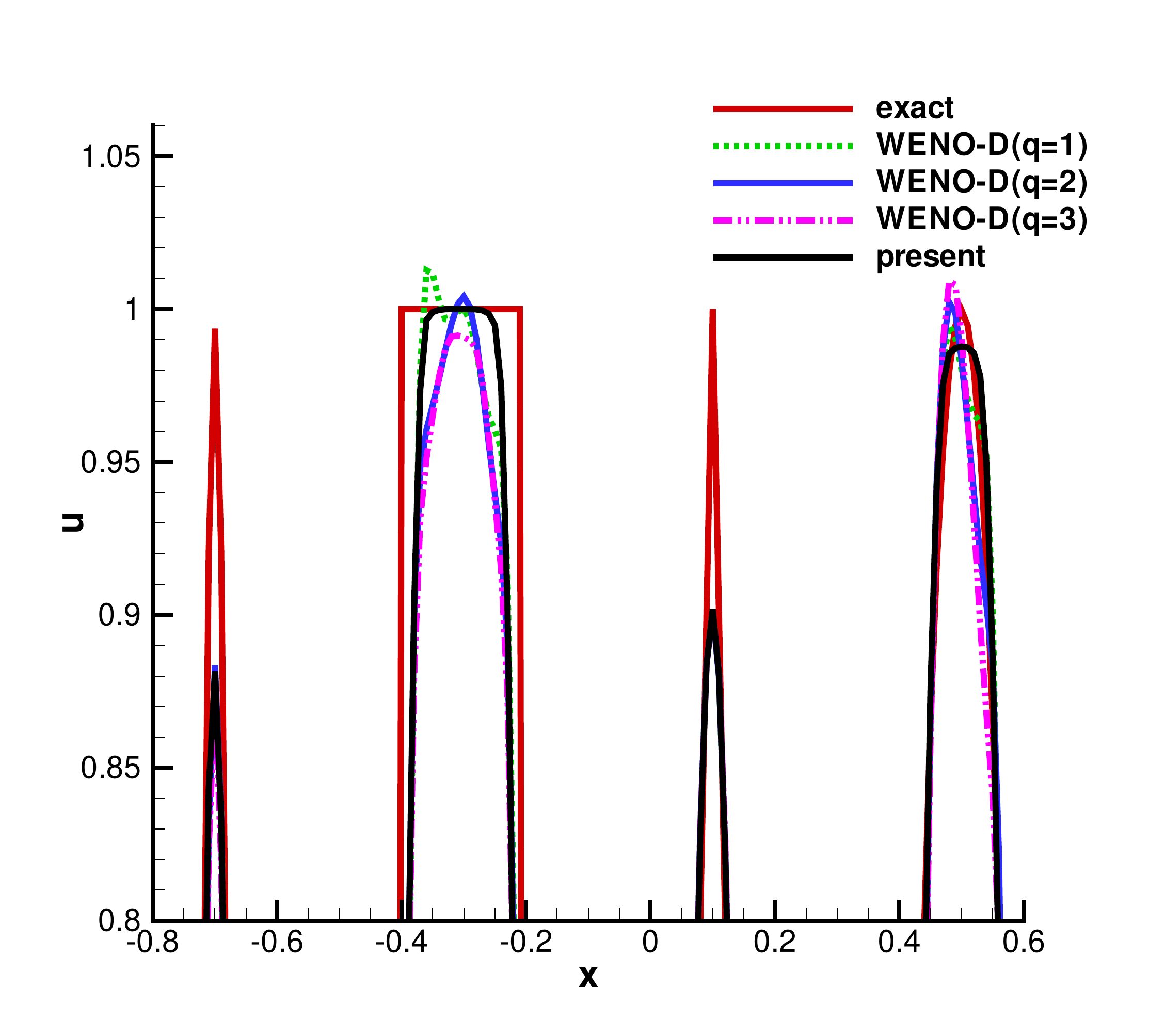}
\caption{Enlarged plot of Fig.\ref{fig:4wave_A1}}
\label{fig:4wave_A100_zoom}
\end{minipage}
\end{figure}

\begin{figure}
\begin{minipage}{0.5\linewidth}
\centering
\includegraphics[width=1.0\textwidth]{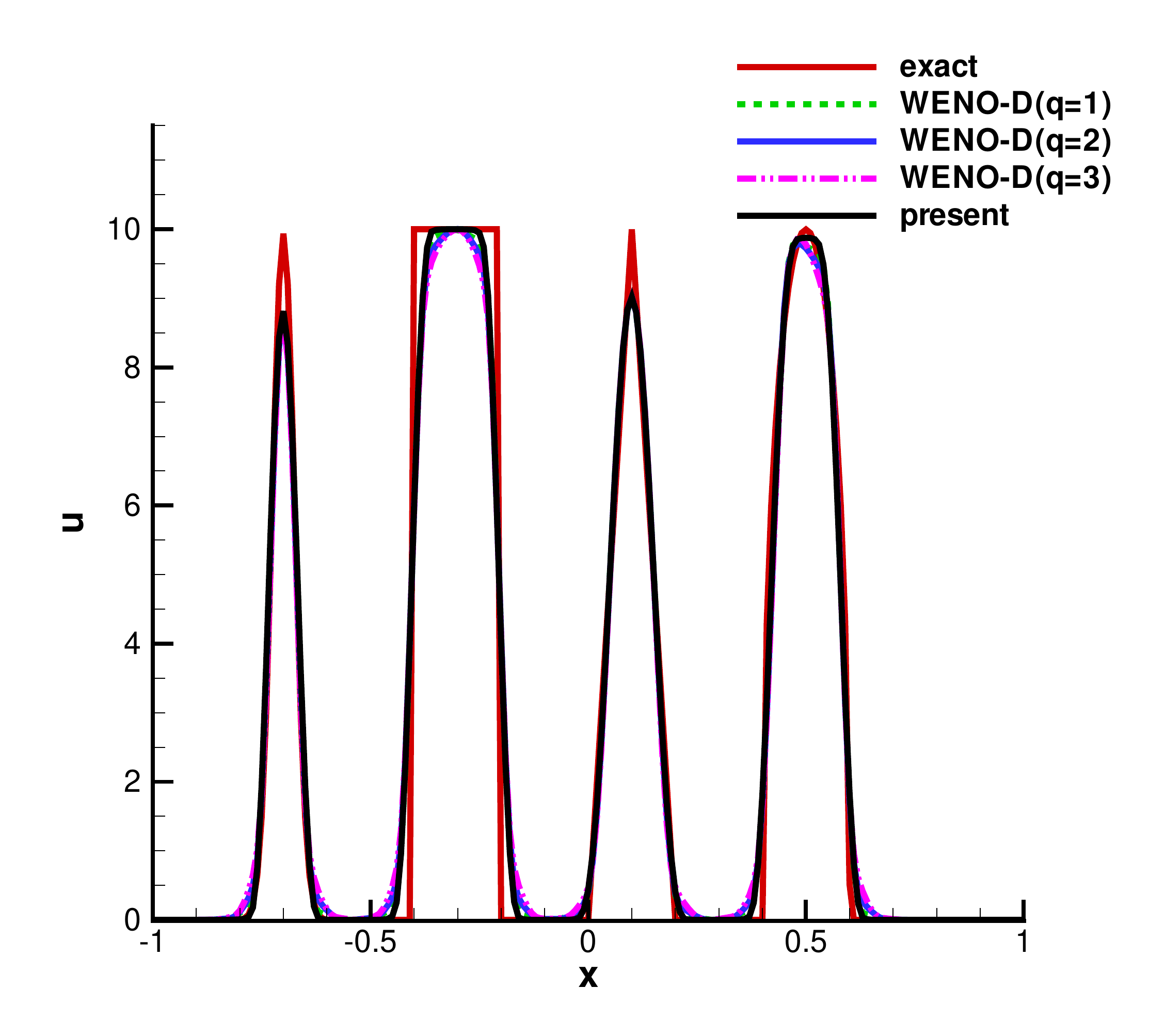}
\caption{Numerical of case 1(\ref{eq:1st_case}), A=0.01}
\label{fig:4wave_A01}
\end{minipage}
\begin{minipage}{0.5\linewidth}
\centering
\includegraphics[width=1.0\textwidth]{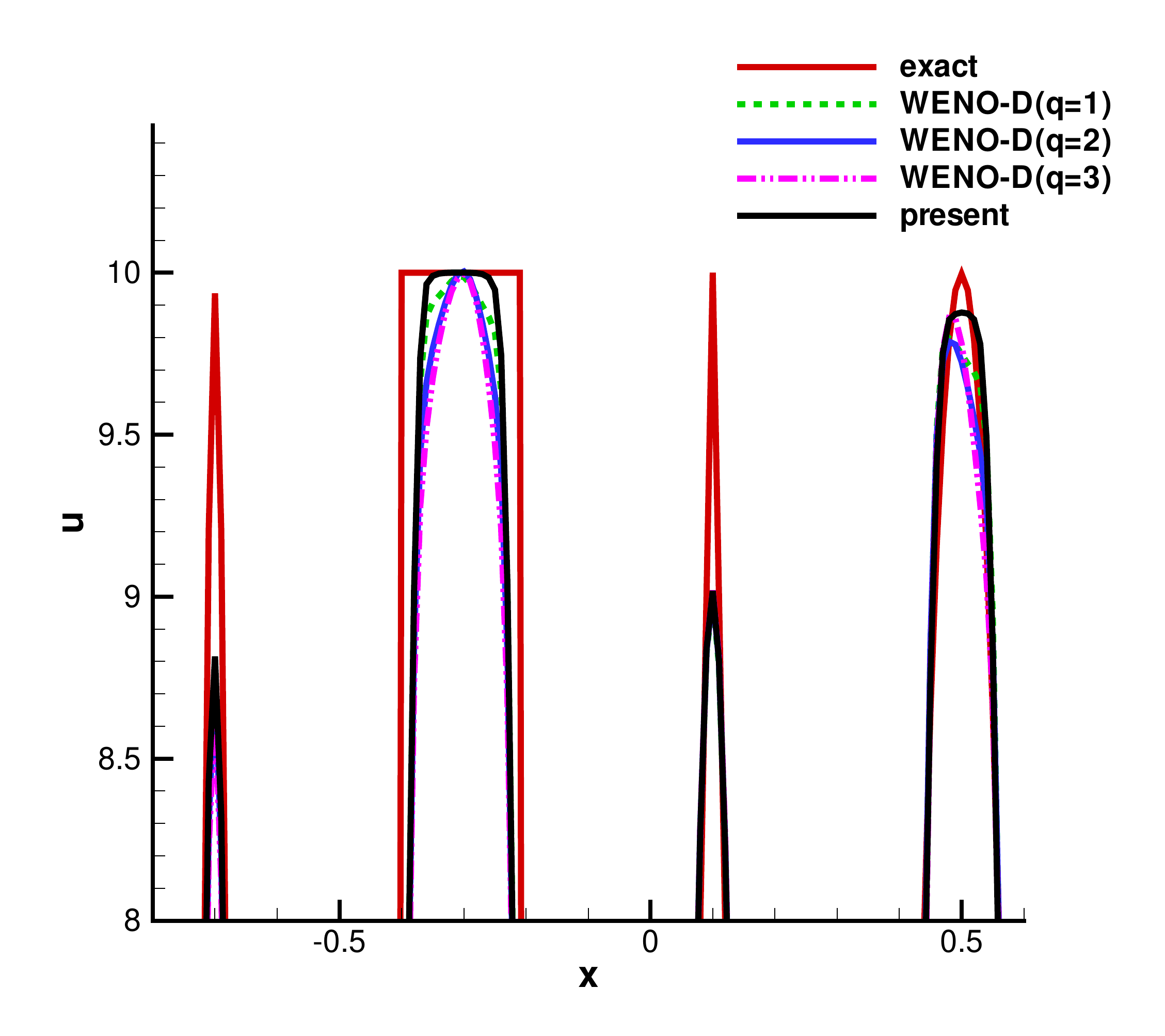}
\caption{Enlarged plot of Fig.\ref{fig:4wave_A01}}
\label{fig:4wave_A01_zoom}
\end{minipage}
\end{figure}

\begin{figure}
\begin{minipage}{0.5\linewidth}
\centering
\includegraphics[width=1.0\textwidth]{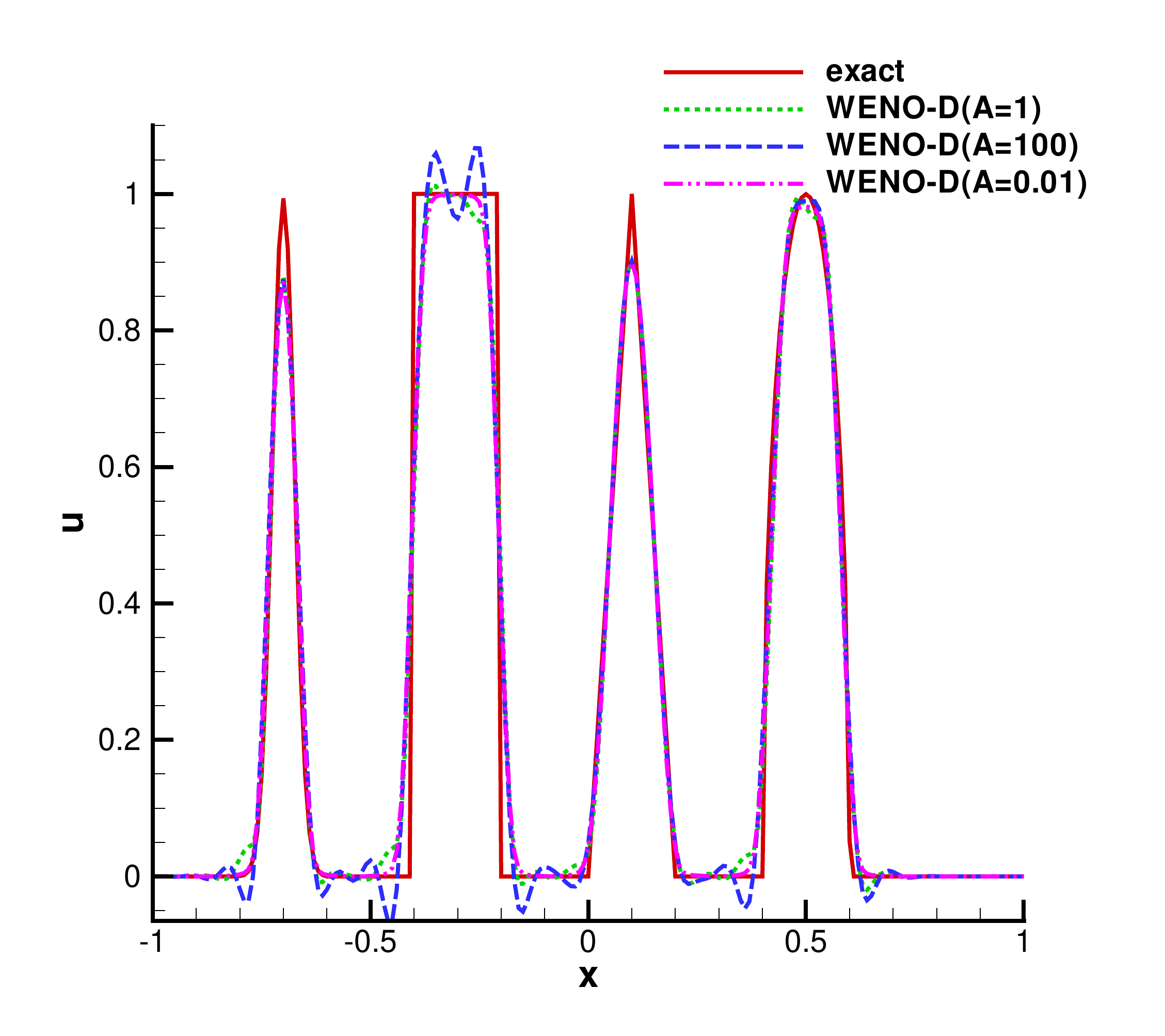}
\caption{Normalized results, WENO-D with $q=1$}
\label{fig:4wave_similarity}
\end{minipage}
\begin{minipage}{0.5\linewidth}
\centering
\includegraphics[width=1.0\textwidth]{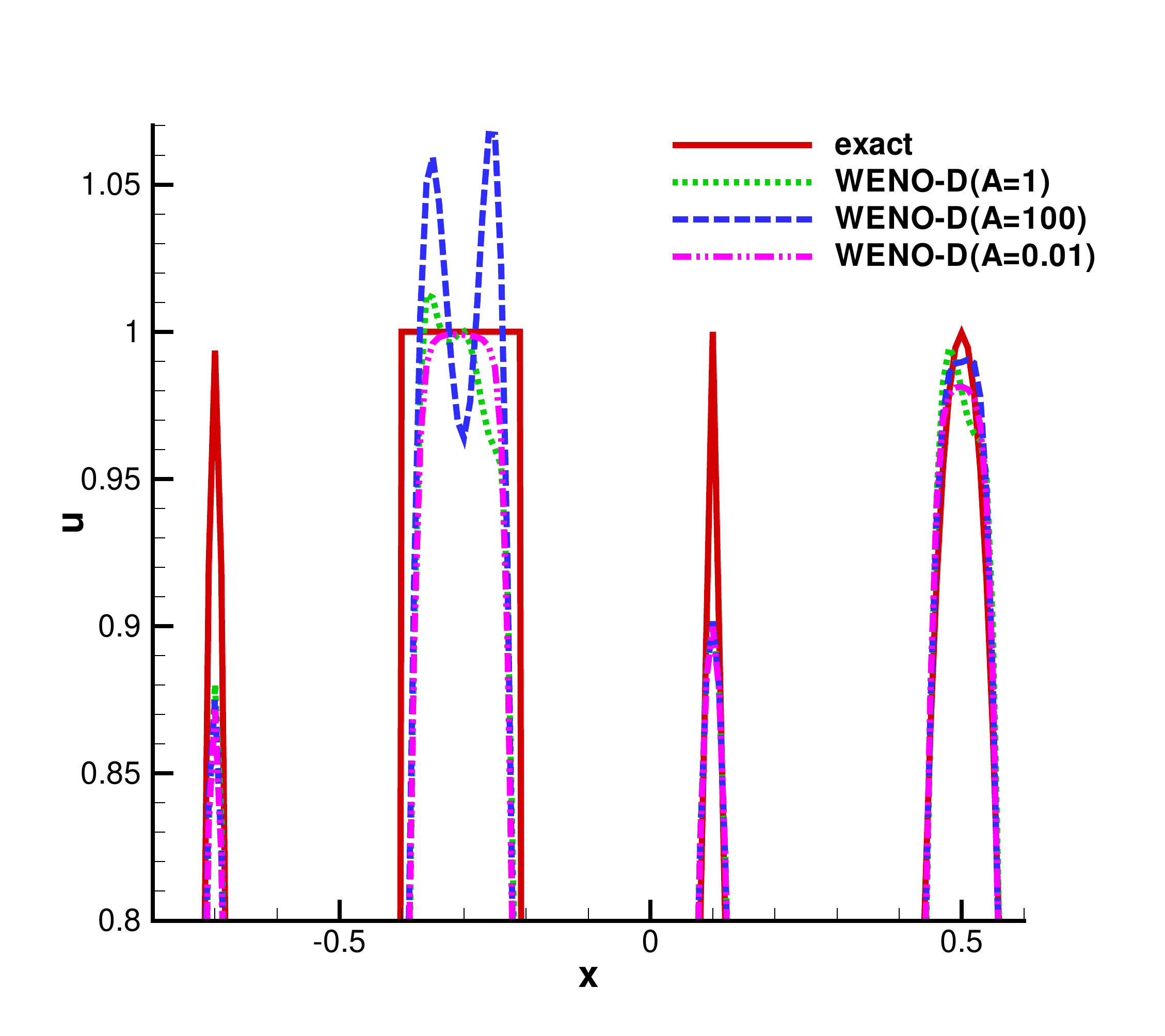}
\caption{Enlarged plot of Fig.\ref{fig:4wave_similarity}}
\label{fig:4wave_similarity_zoom}
\end{minipage}
\end{figure}

{Parameter $A$ in $u_0(x)$ can be regarded as a reference value used to nondimensionalize the unknown variable $u$. In most papers, $A=1$ is always used. In this paper, three cases with $A=100$, $1$ and $0.01$ are tested to demonstrate that how important it is to satisfy the physical requirement that two variables in a addition operation (also in a comparison) should have the same dimension. Numerical results of WENO-D ($\epsilon=10^{-20}$ and $q=1,2,3$) and the present scheme are given in Figs.\ref{fig:4wave_A100}-\ref{fig:4wave_A01_zoom}. It can be seen that, with different values of $A$, the WENO-D scheme shows different spurious behaviours, such as oscillation, apparent asymmetry. As an example, Figs. \ref{fig:4wave_similarity}-\ref{fig:4wave_similarity_zoom} give the comparison of normalized results of WENO-D with $q=1$. It can be seen that, these results lose similarity, even for the smooth ellipse wave. While, the other schemes, including the present scheme, the original WENO-Z scheme\cite{Borges2008} and the WENO-ZA scheme\cite{Liu2018}, keep similarity very well. For compactness, those comparisons are not shown here.}

{Similarly, if a function of the grid spacing is used to replace the parameter $\epsilon$ in the formula of the unnormalized weight $\alpha_k$ of the original WENO-Z scheme\cite{Borges2008}, spurious numerical solutions (such as oscillatory solutions, dissimilar solutions) may also be generated. Such issues caused by unmatched dimensions may emerge when solving the governing equations of fluid dynamics, such as Euler equations or Navier-Stokes equations (Please refer Ref.\cite{Liu2019} for more detailed discussions). Since the purpose of this paper is not to address these issues, applications and comparisons about WENO-D/A are no longer given.}

\begin{figure}
\begin{minipage}{0.5\linewidth}
\centering
\includegraphics[width=1.0\textwidth]{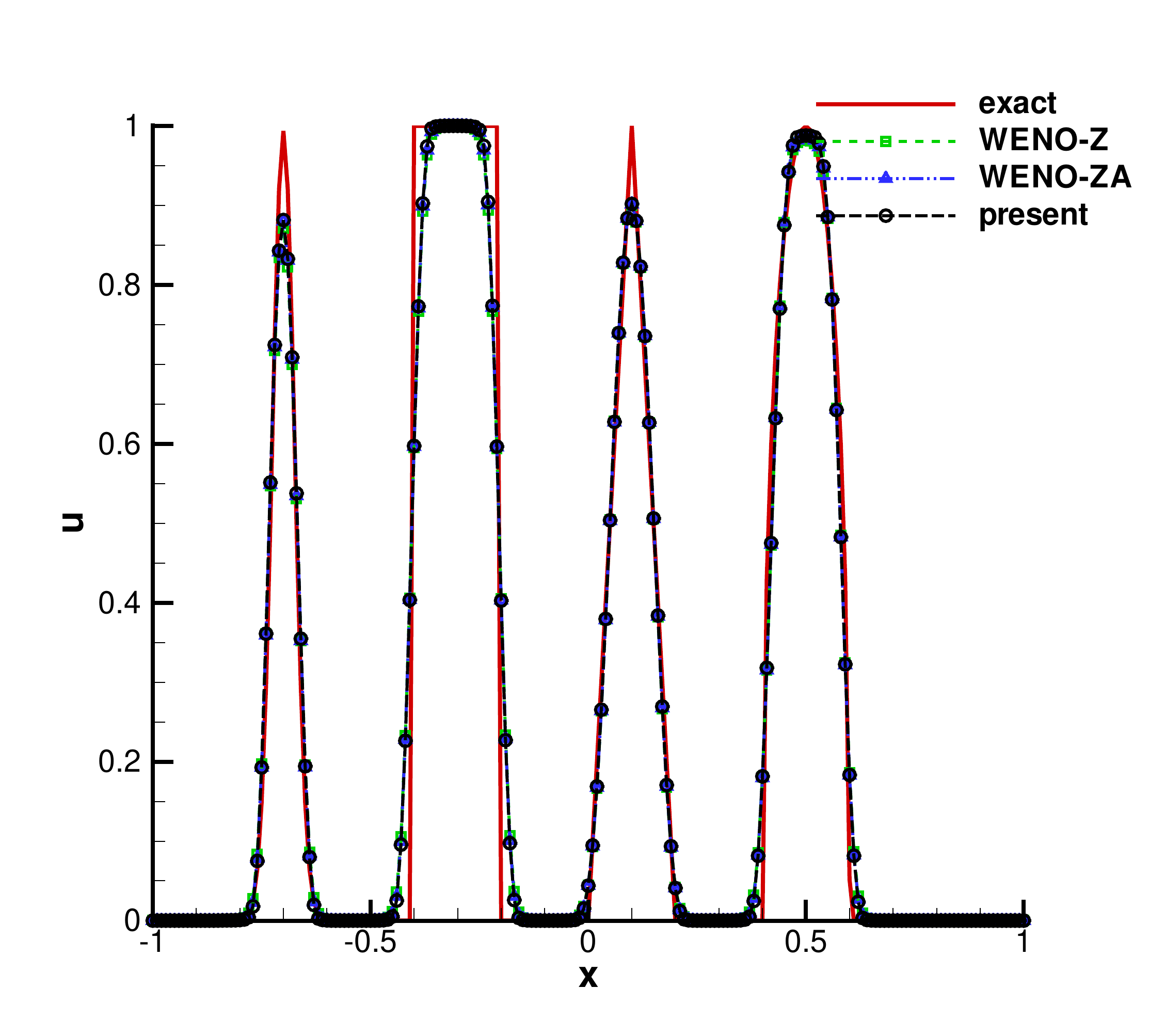}
\caption{Numerical results of case 1(\ref{eq:1st_case}), A=1}
\label{fig:4wave}
\end{minipage}
\begin{minipage}{0.5\linewidth}
\centering
\includegraphics[width=1.0\textwidth]{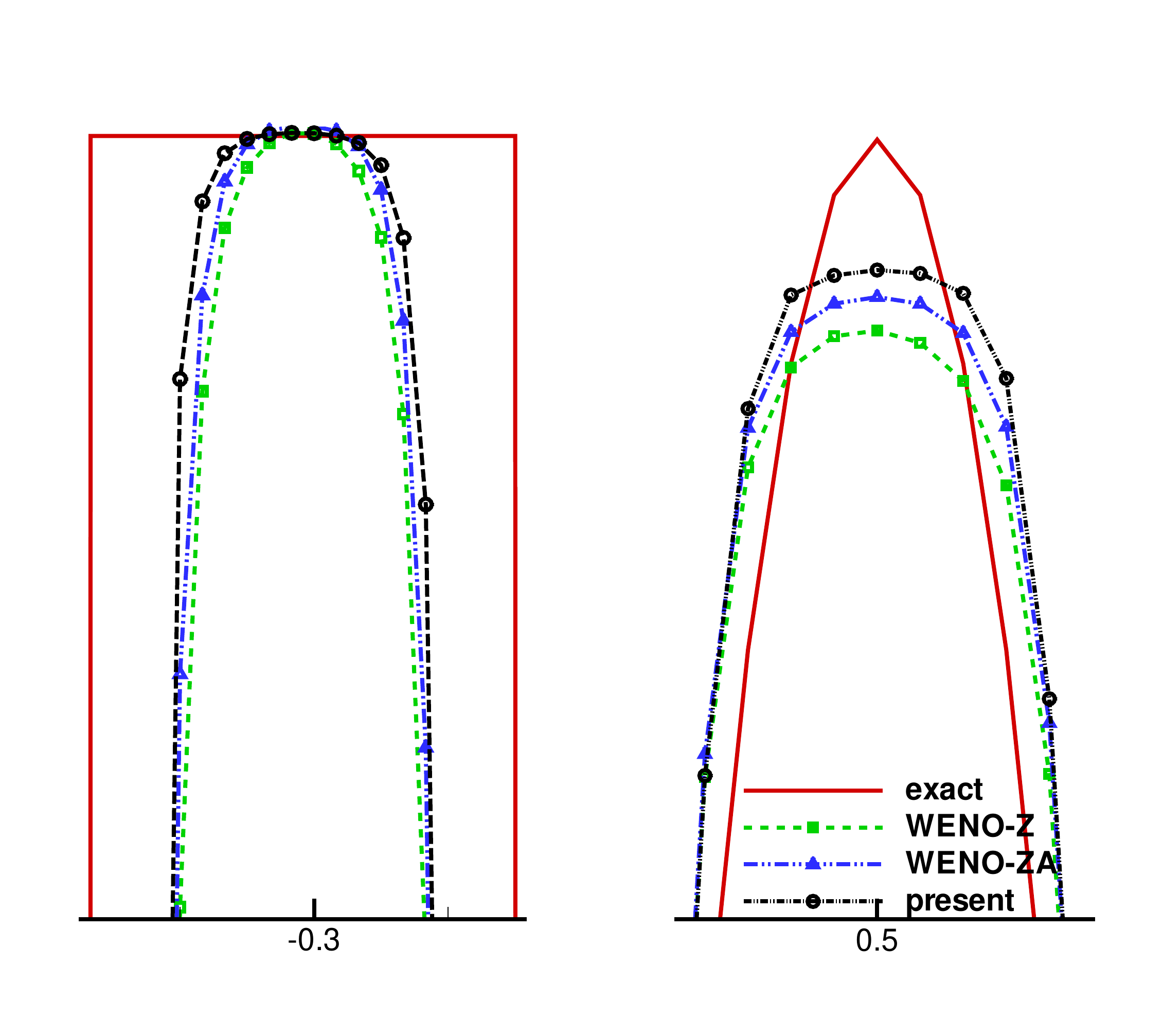}
\caption{Enlarged plot of Fig.\ref{fig:4wave}}
\label{fig:4wave11}
\end{minipage}
\end{figure}

Figs.\ref{fig:4wave} and \ref{fig:4wave11} are the numerical comparisons of the present scheme, the WENO-Z scheme and the WENO-ZA scheme. It can be seen that, the present scheme resolves both discontinuity(the square wave) and smooth solution (ellipse wave) more accurate than WENO-Z and WENO-ZA.

Case 2: The initial condition is given in \eqref{eq:single}.
$$
u_0(x)=\left\{
\begin{array}{ll}
-\text{sin}(\pi x)-\dfrac{1}{2}x^3,&-1\le x< 0, \\
-\text{sin}(\pi x)-\dfrac{1}{2}x^3+1,&0\le x \le 1.
\end{array}\right .
$$

\begin{figure}
\begin{minipage}{0.5\linewidth}
\centering
\includegraphics[width=1.0\textwidth]{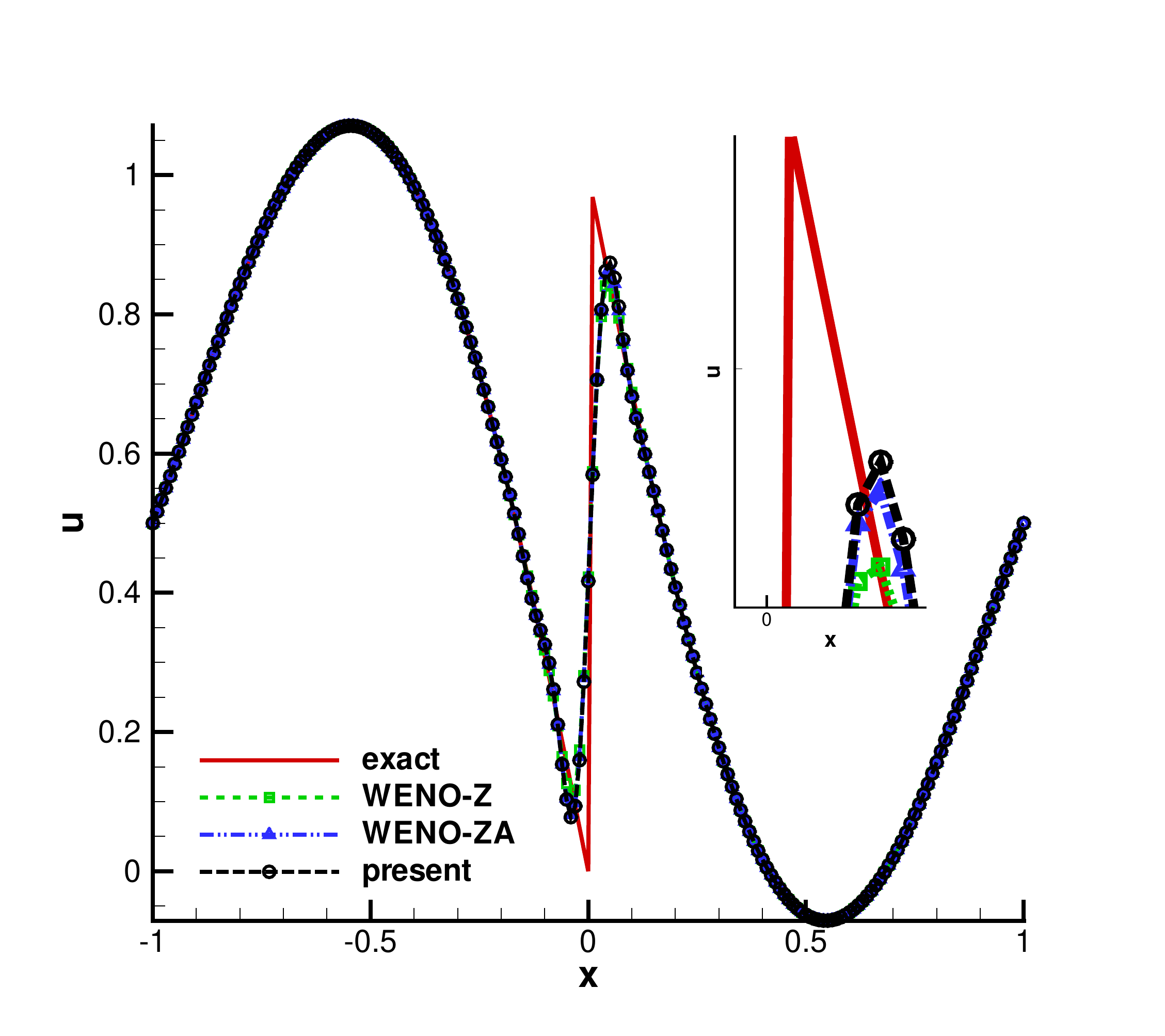}
\caption{Numerical results of case 2(\ref{eq:single}) }
\label{fig:1wave}
\end{minipage}
\begin{minipage}{0.5\linewidth}
\centering
\includegraphics[width=1.0\textwidth]{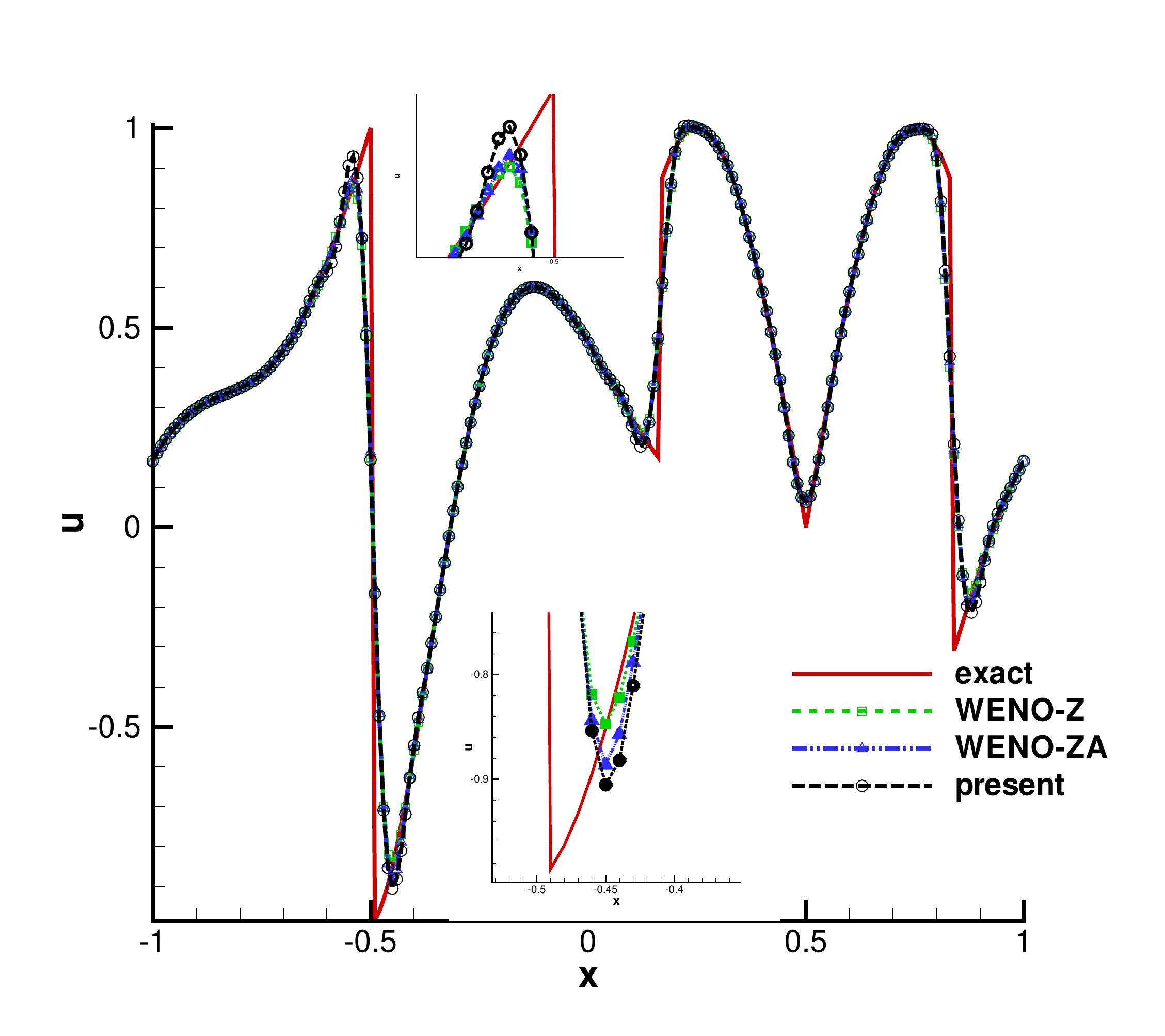}
\caption{Numerical results of case 3(\ref{eq:3rd_case})}
\label{fig:2wave}
\end{minipage}
\end{figure}

Fig.\ref{fig:1wave} shows the numerical results at $T=6$ with $N=200$. The new method improves the resolution near discontinuity.

Case 3: The initial condition is given as
\begin{align}\label{eq:3rd_case}
u_0(x)=\left\{
\begin{array}{ll}
-xsin(3\pi x^2/2), & \ -1\le x<-1/3\\
|sin(2\pi x)|, & \ -1/3\le x\le1/3\\
2x-1-sin(3\pi x)/6, & \ otherwise
\end{array}
\right .
\end{align}

The results at $T=6$ with $N=200$ are plotted in Fig.\ref{fig:2wave}. Near discontinuity, the present scheme performs better than the other schemes.

\subsection{One-dimensional Euler problems}

The governing equations are as follows
\begin{equation}
\frac{\partial{U}}{\partial{t}}+\frac{\partial{F}}{\partial{x}}=0,
\end{equation}
where $U=(\rho,\rho{u},E)^T,\ F(U)=\big(\rho{u},\rho{u^2}+p,u(E+p)\big)^T$, $\rho,u,E\  \text{and} \  p$ are the density, the velocity, the total energy, and the pressure respectively. The equation of state is given by $E=\dfrac{p}{\gamma-1}+\dfrac{1}{2}\rho{u^2}$, where $\gamma=1.4$ is the ratio of specific heat. Time step is taken as
\begin{equation}
\Delta{t}=\frac{\sigma\Delta{x}}{max_i\left(\left|u_i\right|+c_i\right)},
\end{equation}
where, $\sigma=0.5$ is the CFL number. $c$ is the speed of sound and given by $c=\sqrt{\gamma{p}/\rho}$. The LF flux-splitting method is used and the WENO reconstruction is carried out in local characteristic fields \cite{2009esweno}.
The Shu-Osher problem \cite{Shu1988} and the interactive blast waves problem \cite{Borges2008} are calculated. All the reference solutions are obtained by the WENO-Z scheme with a grid of 2000.

\subsubsection{Shu-Osher problem}
\begin{figure}
\begin{minipage}{0.5\linewidth}
\centering
\includegraphics[width=1.0\textwidth]{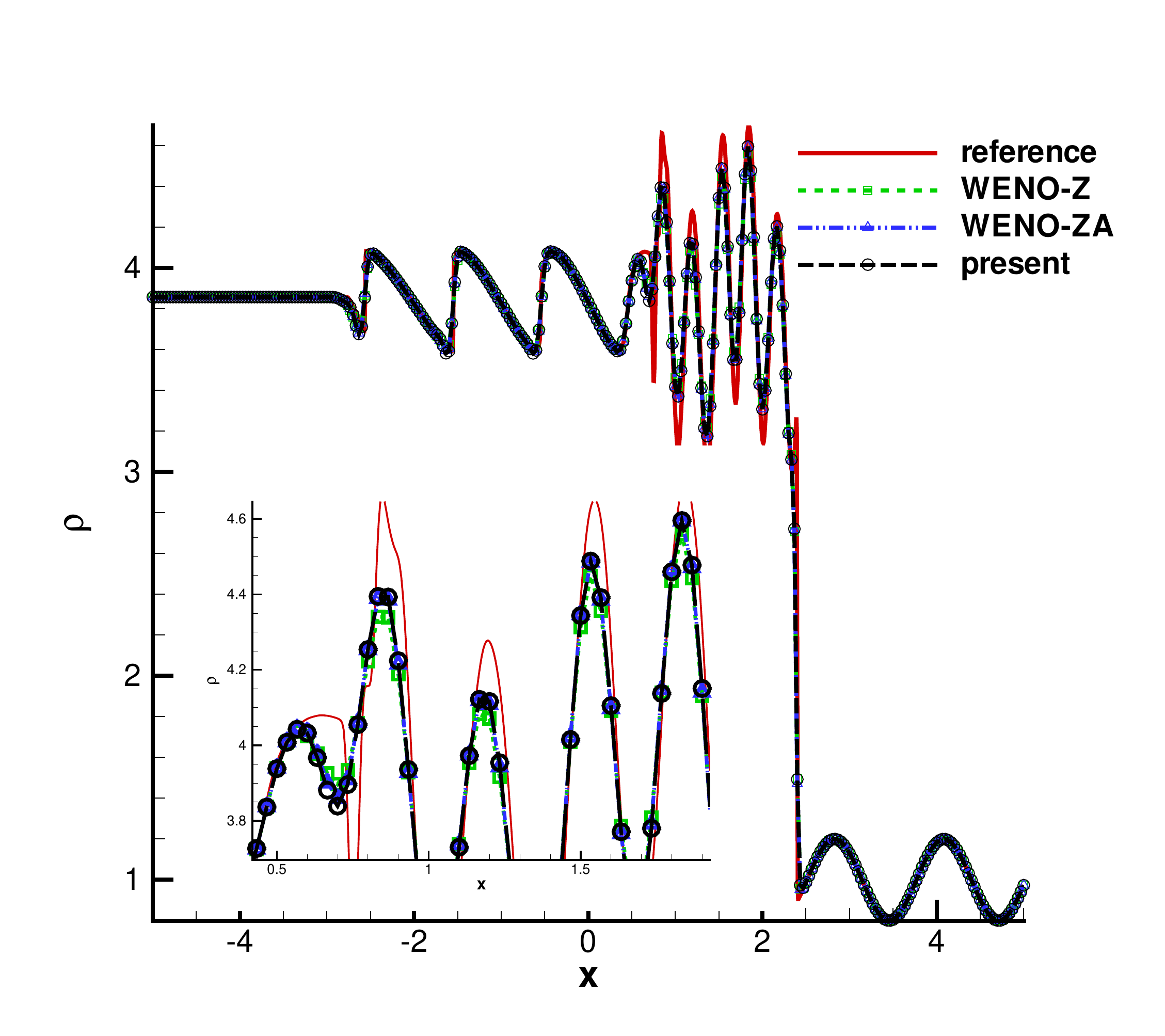}
\caption{Density distribution of Case 4, Shu-Osher problem}
\label{fig:shu}
\end{minipage}
\begin{minipage}{0.5\linewidth}
\centering
\includegraphics[width=1.0\textwidth]{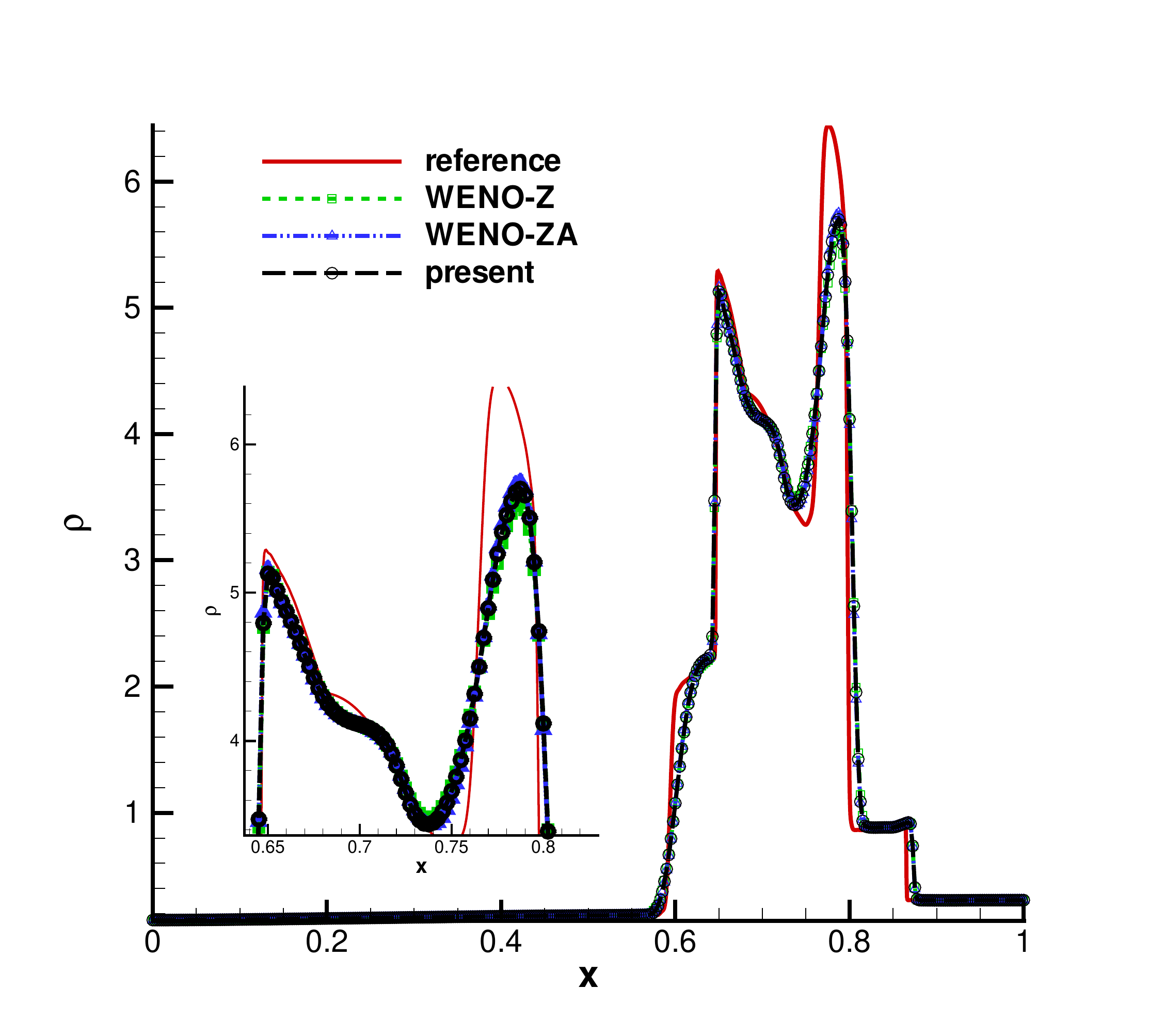}
\caption{Density distribution of Case 5, interaction of two blast waves}
\label{fig:2-blast}
\end{minipage}
\end{figure}

The first 1-D case is the Shu-Osher problem \cite{Shu1988} with the initial condition
\begin{align}\label{eq:shu}
(\rho,u,p)=
\begin{cases}
(3.857143,\ 2.629369,\ 31/3) & -5\le{x}<-4, \\
(1+0.2\text{sin}(5x),\ 0,\ 1) & -4\le{x}\le{5} .
\end{cases}&
\end{align}
Fig.\ref{fig:shu} gives the distributions of density at $t=1.8$ by using $N=300$. As this figure shows, the present scheme is almost the same as the WENO-ZA scheme, both of
them are less dissipative than WENO-Z.

\subsubsection{Two interacting blast waves}
The second 1D-case is the interactive blast waves problem \cite{Borges2008} with the initial condition
\begin{align}
(\rho,u,p)=
\begin{cases}
(1,\ 0,\ 1000) & 0\le{x}<0.1, \\
(1,\ 0,\ 0.01) & 0.1\le{x}<0.9, \\
(1,\ 0,\ 100) & 0.9\le{x}\le{1} .
\end{cases}&
\end{align}
The numerical results at $t=0.038$ with $N=400$ are presented in Fig.\ref{fig:2-blast}. The three schemes can capture strong shocks well.

\subsection{Two-dimensional Euler problems}

The governing equation is the two-dimensional Euler equations
\begin{equation}\label{eq:2d}
\frac{\partial{U}}{\partial{t}}+\frac{\partial{F}}{\partial{x}}+\frac{\partial{G}}{\partial{y}}=0,
\end{equation}
where the conservative variables $U$ and the inviscid flux vectors $F$ and $G$ are
\begin{equation}
\begin{split}
&U=\left[\begin{matrix}
\rho\\
\rho u\\
\rho v\\
E\\
\end{matrix}
\right],\
F=\left[\begin{matrix}
\rho u\\
\rho u^2+p\\
\rho uv\\
Eu+pu\\
\end{matrix}
\right],\
G=\left[\begin{matrix}
\rho v\\
\rho uv\\
\rho v^2+p\\
Ev+pv\\
\end{matrix}
\right].
\end{split}
\end{equation}
The energy is given by
\begin{equation}
E=\frac{p}{\gamma-1}+\frac{\rho}{2}(u^2+v^2).
\end{equation}
The Steger-Warming flux vector splitting method \cite{Steger1981} is used for the inviscid convective fluxes, and the time step is taken as follows
\begin{equation}
\begin{split}
&\Delta{t}=\sigma\frac{\Delta{t_x}\Delta{t_y}}{\Delta{t_x}+\Delta{t_y}}, \\
&\Delta{t_{x}}=\dfrac{\Delta{x}}{max_{i,j}\left(\left|u_{i,j}\right|+c_{i,j}\right)}, \\
&\Delta{t_y}=\dfrac{\Delta{y}}{max_{i,j}\left(\left|v_{i,j}\right|+c_{i,j}\right)}.
\end{split}
\end{equation}

\subsubsection{Riemann problems}

Two-dimensional Riemann problems with different initial configurations have been extensively employed to examine the numerical schemes for Euler equations\cite{Schulz_Rinne1993, Kurganov2002, Abedian2014, Jung2017, Deng2019}. Two cases are calculated in this section.

\begin{figure}
\centering
\includegraphics[width=0.8\textwidth]{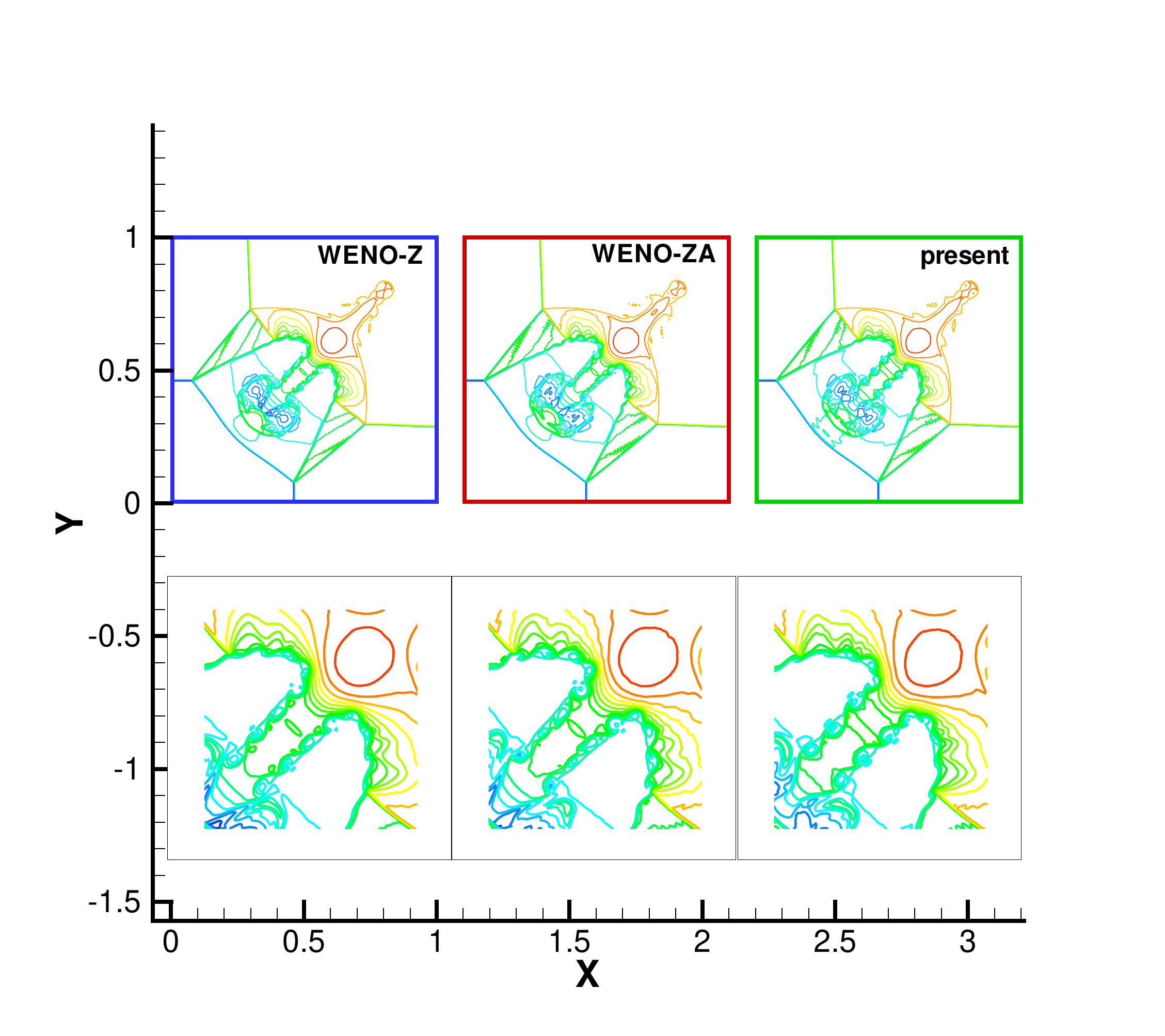}
\caption{Density contours for the first 2D Riemann problem \eqref{eq:RM1}, $400\times400$}
\label{fig:RM4}
\end{figure}

Case 1: The initial conditions of the first case are given as
\begin{align}\label{eq:RM1}
(\rho,u,v,p)=
\begin{cases}
(1.5,\ 0,\ 0,\ 1.5) & 0.8\le{x}\le{1},0.8\le{y}\le{1}, \\
(0.5323,\ 1.206,\ 0,\ 0.3) & 0\le{x}<0.8,0.8\le{y}\le{1}, \\
(0.138,\ 1.206,\ 1.206,\ 0.029) & 0\le{x}<0.8,0\le{y}<0.8, \\
(0.5323,\ 0,\ 1.206,\ 0.3) & 0.8\le{x}\le{1},0\le{y}<0.8.
\end{cases}&
\end{align}
A grid of $400\times400$ is used. The density contours at $t=0.8$ are shown in Fig.\ref{fig:RM4}. It can be seen that the three schemes can capture
reflection shocks and contact discontinuities well.
But the present scheme can resolve the roll-ups of the Kelvin{-}Helmholtz instability with finer structures than the other two schemes.

\begin{figure}
\centering
\subfigure[WENO-Z]{
\includegraphics[width=0.45\textwidth]{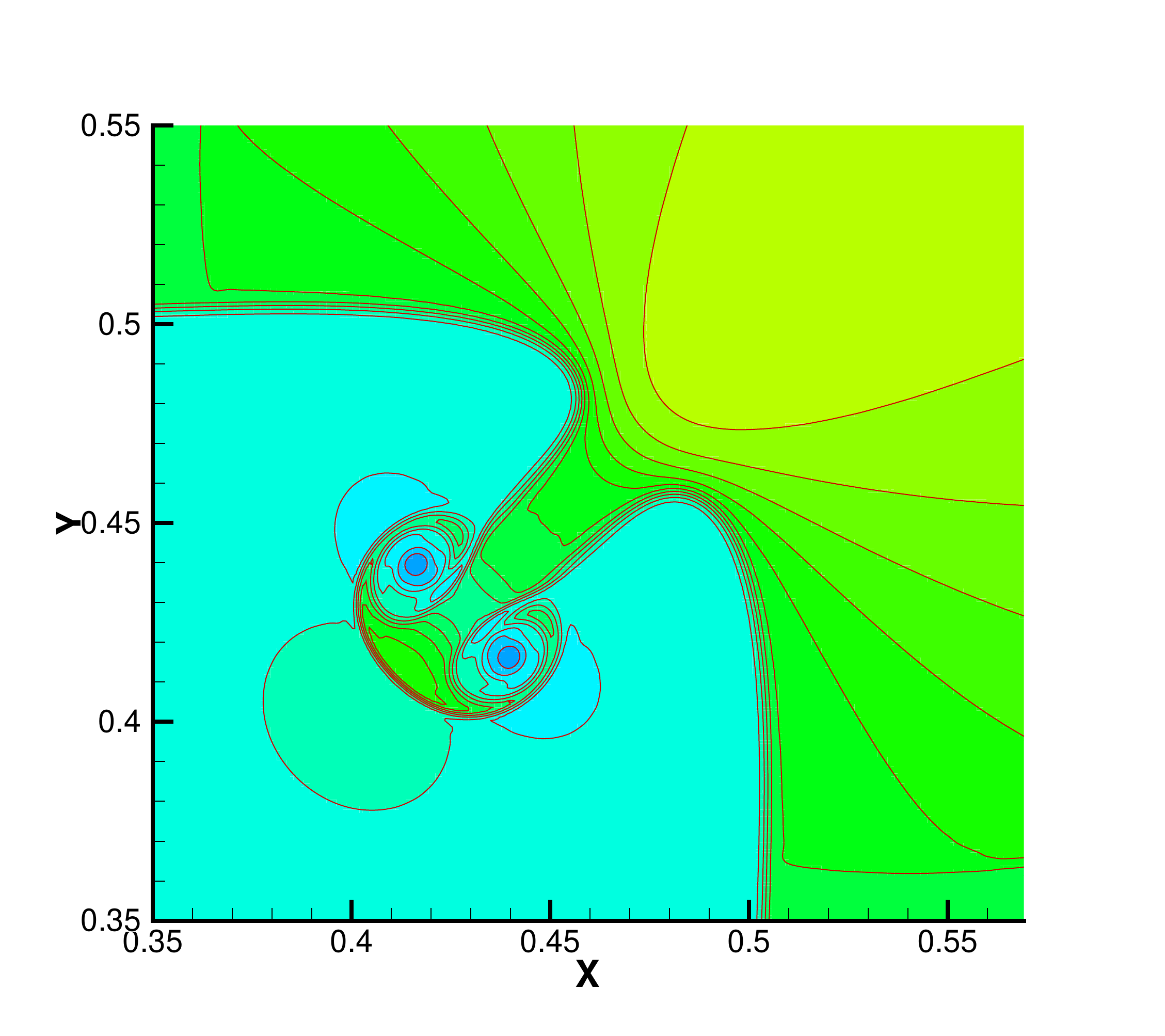}
}
\subfigure[WENO-ZA]{
\includegraphics[width=0.45\textwidth]{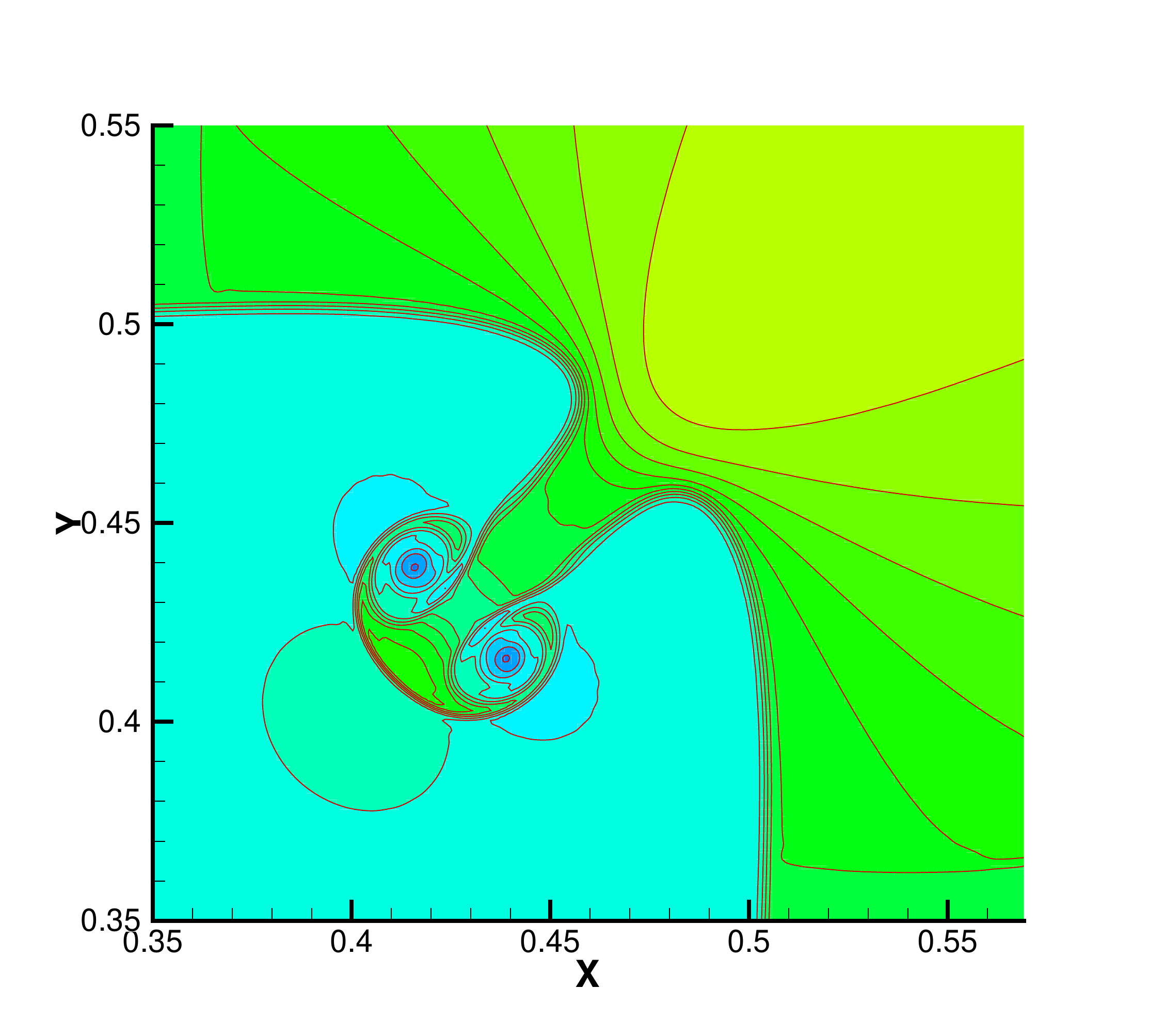}
}
\subfigure[present]{
\includegraphics[width=0.45\textwidth]{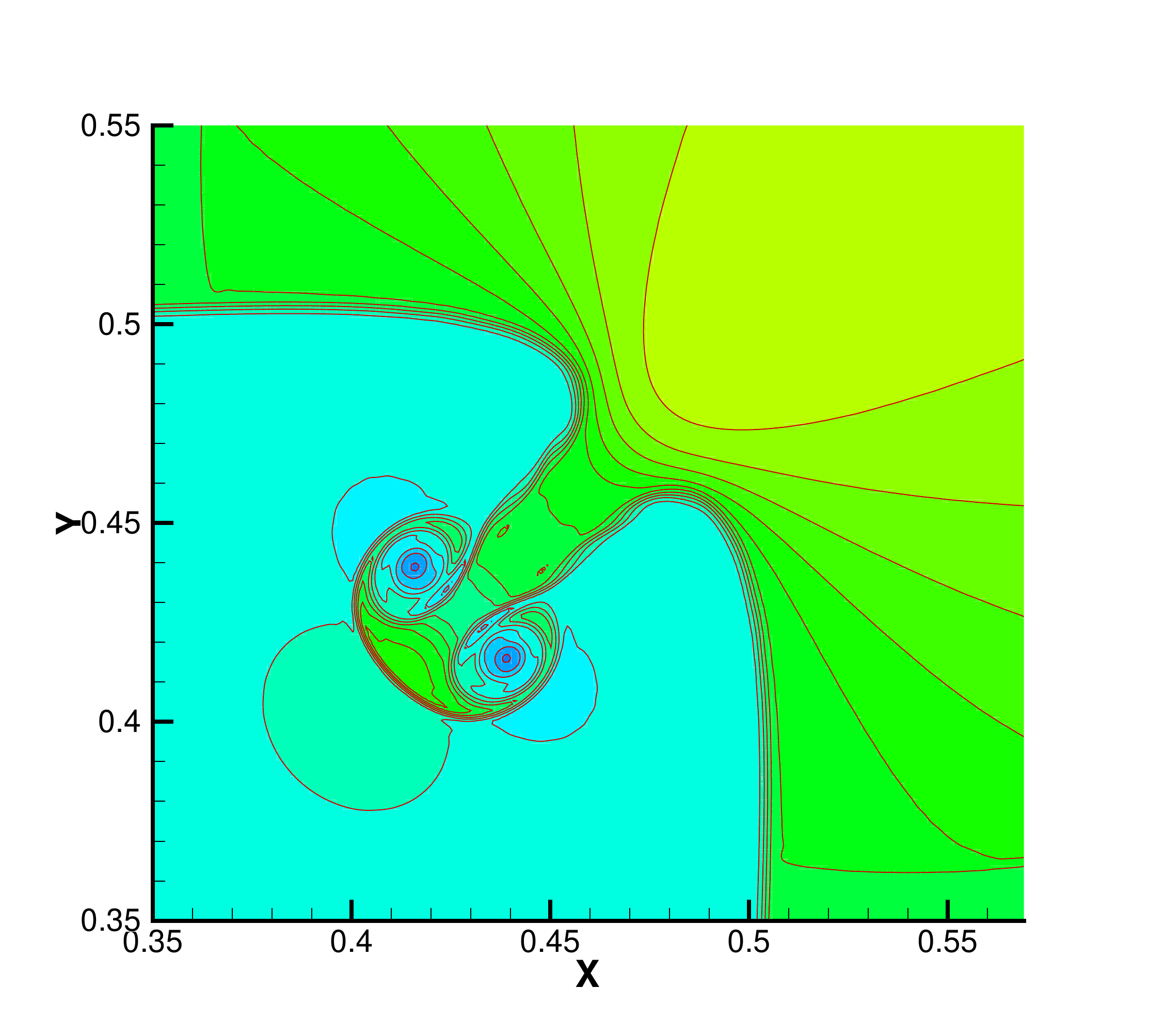}
}
\caption{Density contours for the second 2D Riemann problem \eqref{eq:RM2}. Grid: $1200\times 1200$.}
\label{fig:RM1200}
\end{figure}

\begin{figure}
  \centering
\subfigure[WENO-Z]{
\includegraphics[width=0.45\textwidth]{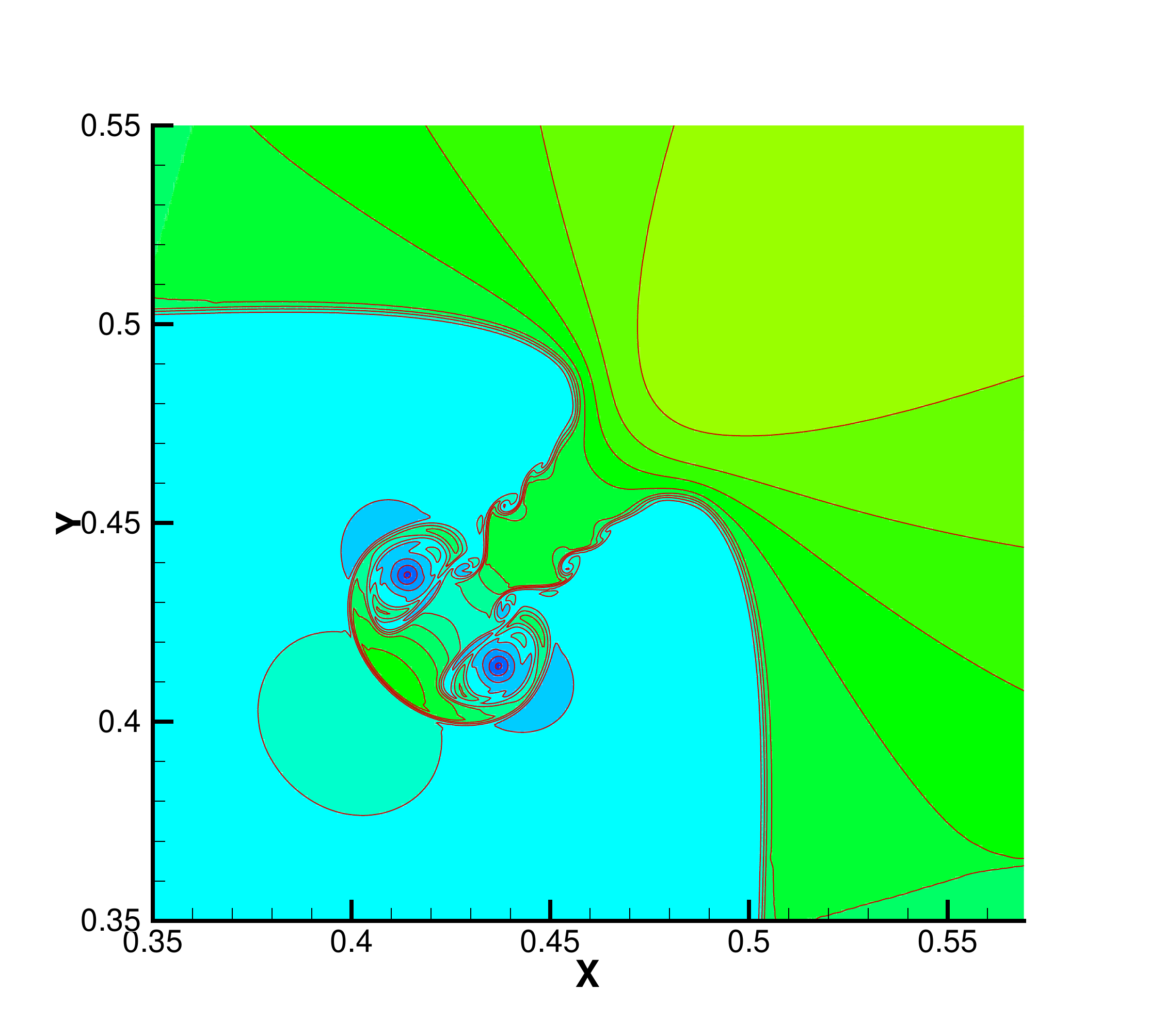}
}
\subfigure[WENO-ZA]{
\includegraphics[width=0.45\textwidth]{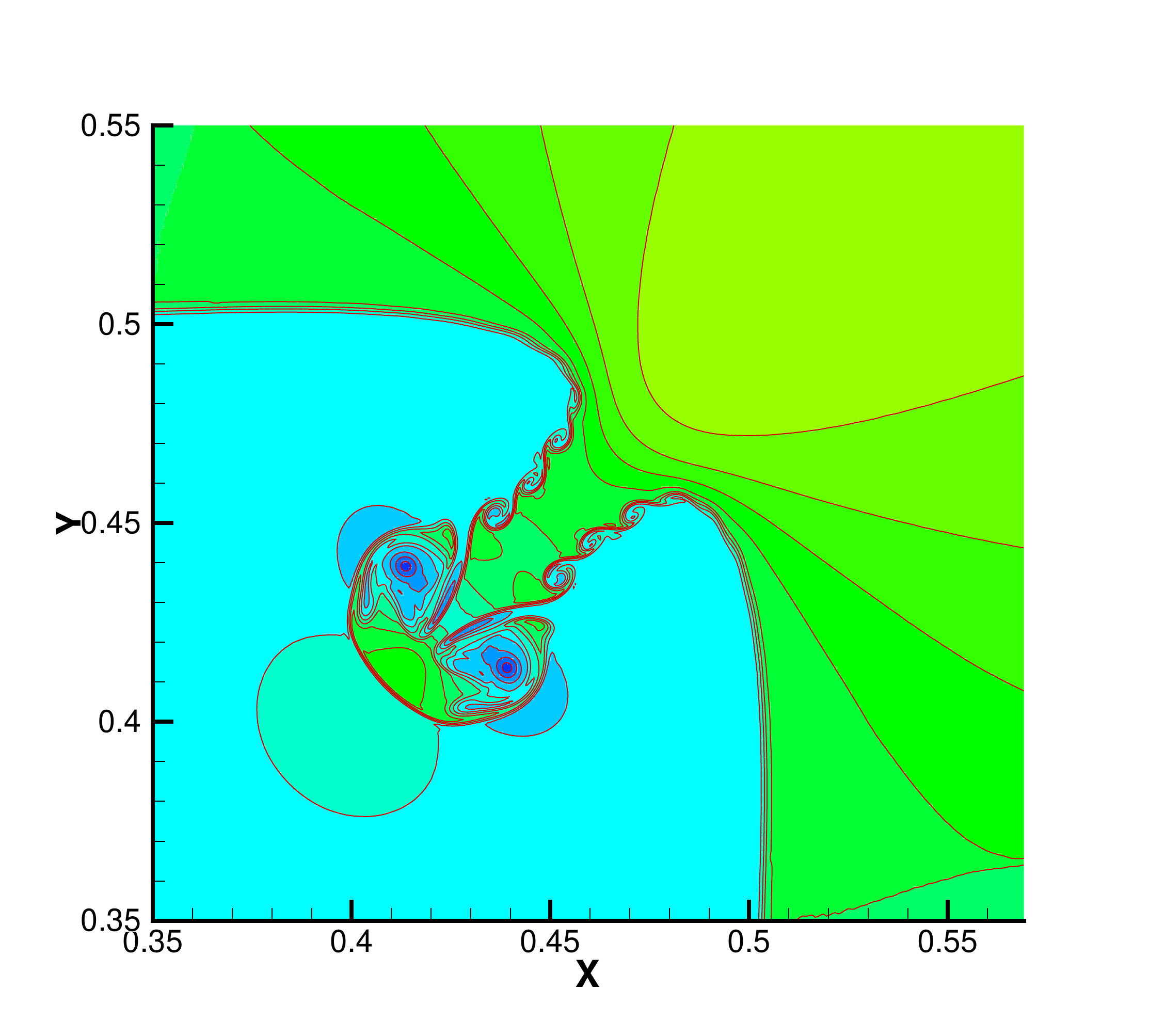}
}
\subfigure[present]{
\includegraphics[width=0.45\textwidth]{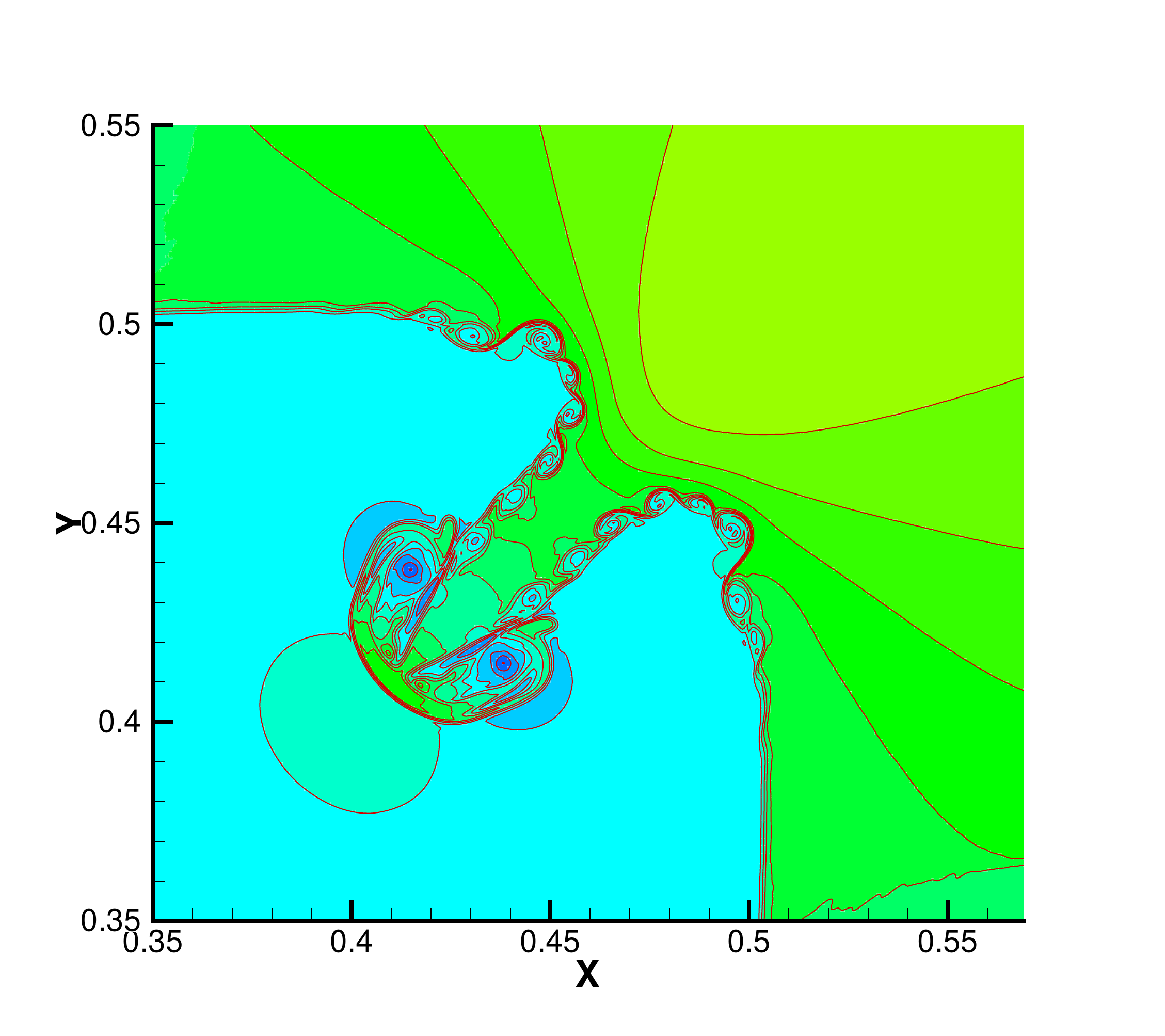}
}
\caption{Density contours for the second 2D Riemann problem \eqref{eq:RM2}. Grid: $2400\times 2400$.}
\label{fig:RM2400}
\end{figure}

Case 2: The initial conditions of the second case are
\begin{align}\label{eq:RM2}
(\rho,u,v,p)=
\begin{cases}
(0.8,\ 0,\ 0,\ 1.0) & {x}\le{0.5},{y}\le{0.5}, \\
(1.0,\ 0.7276,\ 0,\ 1.0) & {x}\le0.5,{y}>{0.5}, \\
(1.0,\ 0.0,\ 0.7276,\ 1.0) & {x}>0.5,{y}\le 0.5, \\
(0.5313,\ 0.0,\ 0.0,\ 0.4) & {x}>0.5,{y}>0.5.
\end{cases}&
\end{align}
For this case, the fine structures of the KH instability along the slip lines are hardly reproduced unless high-order
schemes with minimized numerical dissipation or very fine computational grids are used\cite{Deng2019}. In our calculations, two sets of grids of
$1200\times 1200$ and $2400\times 2400$ are tested. Density contours at $t=0.25$ are plotted in Fig.\ref{fig:RM1200} and \ref{fig:RM2400}. With grid of $1200\times 1200$, three schemes almost can not resolve the small-scale structures generated by the KH instability. With the finer grid, the small-scale structures are generated, and their richness indicates the presented scheme has the lowest numerical dissipation.

\subsubsection{Rayleigh-Taylor instability}

\begin{figure}
\centering
\includegraphics[width=0.7\textwidth]{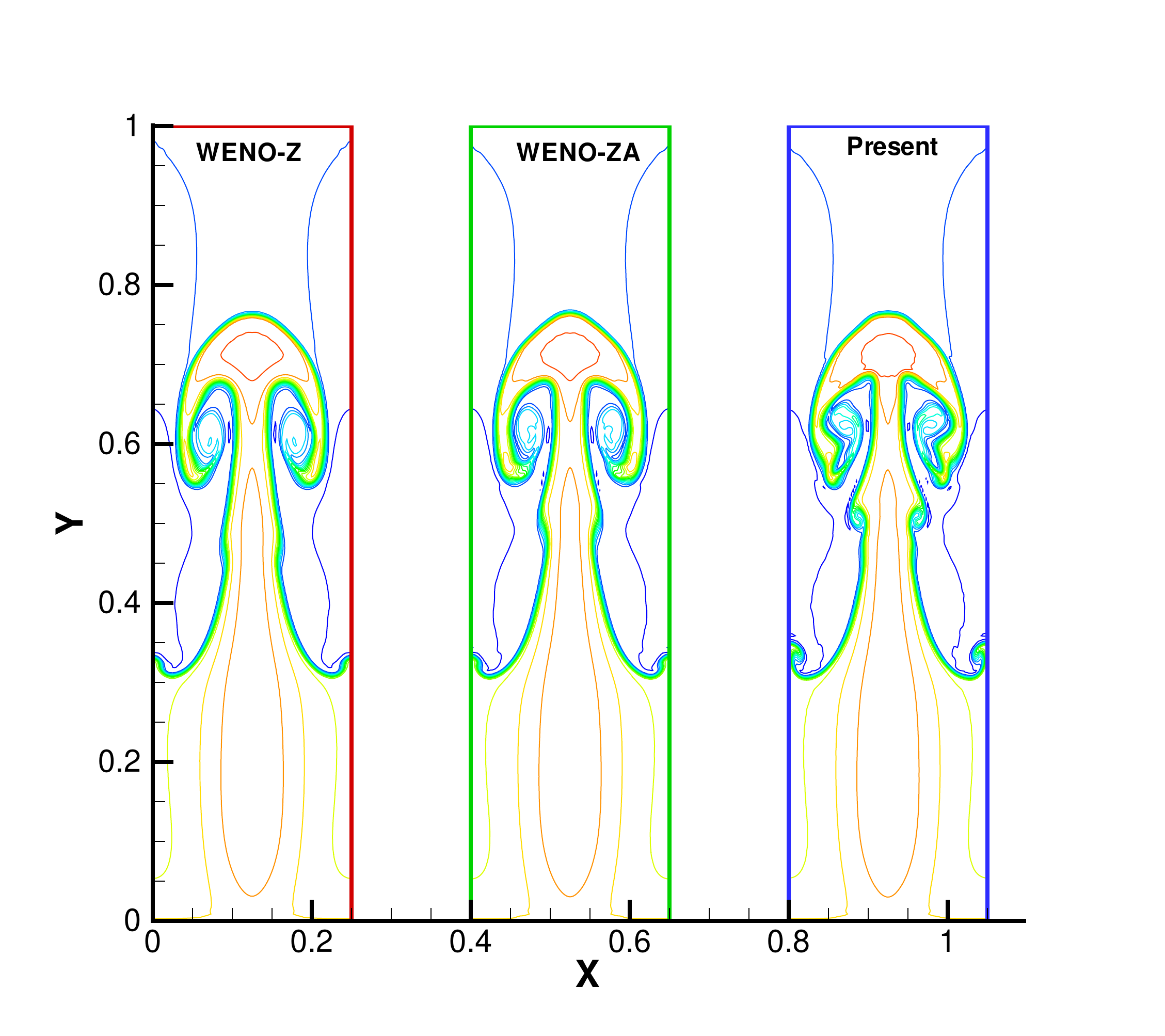}
\caption{Density contours of Rayleigh-Taylor instability, $120\times480$}
\label{fig:RT8}
\end{figure}

The two-dimensional Rayleigh-Taylor instability problem \cite{Shi2003, Yong2001} is often used to assess the dissipation property of a high-order scheme. It describes the interface instability between fluids with different densities when acceleration is directed from a heavy fluid to a light one. The gravitational effect is introduced by adding $\rho$ and $\rho{v}$ to the flux of the $y$-momentum and the energy equations, respectively. The initial distribution is
\begin{align}\label{eq:40}
&(\rho,u,v,p)=
\begin{cases}
(2,0,-0.025\alpha{\text{cos}(8\pi{x})},2y+1),& 0\le{y}<1/2, \\
(1,0,-0.025\alpha{\text{cos}(8\pi{x})},y+3/2),& 1/2\le{y}<1,
\end{cases}&
\end{align}
and $\alpha=\sqrt{\gamma{p}/\rho}$ is the speed of sound with $\gamma=5/3$. The computational domain is $[0,0.25]\times{[0,1]}$. The left and right boundaries are reflective boundary conditions, and the top and bottom boundaries are set as $(\rho,u,v,p)=(1,0,0,2.5)$ and $(\rho,u,v,p)=(2,0,0,1)$, respectively.
The solution at $t=1.95$ is solved with a mesh of $120\times 480$.
The density contours are plotted in Fig.\ref{fig:RT8}.
As observed in previous cases, due to lower dissipation, the present scheme resolves more clear unstable structures than the other two schemes.

\subsubsection{Double Mach reflection}

The double Mach reflection problem describes the reflection of a planar Mach shock in air hitting a
wedge \cite{Woodward1984}. The initial conditions are given as
\begin{align}\label{eq:DM}
&(\rho,u,v,p)=
\begin{cases}
(8,0,8.25cos(\pi/6),-8.25sin(\pi/6),116.5),& x<\frac{1}{6}+\frac{y}{\sqrt{3}}, \\
(1.4,0,0,1.0),& x<\frac{1}{6}+\frac{y}{\sqrt{3}}.
\end{cases}&
\end{align}

\begin{figure}
\centering
\includegraphics[width=0.8\textwidth]{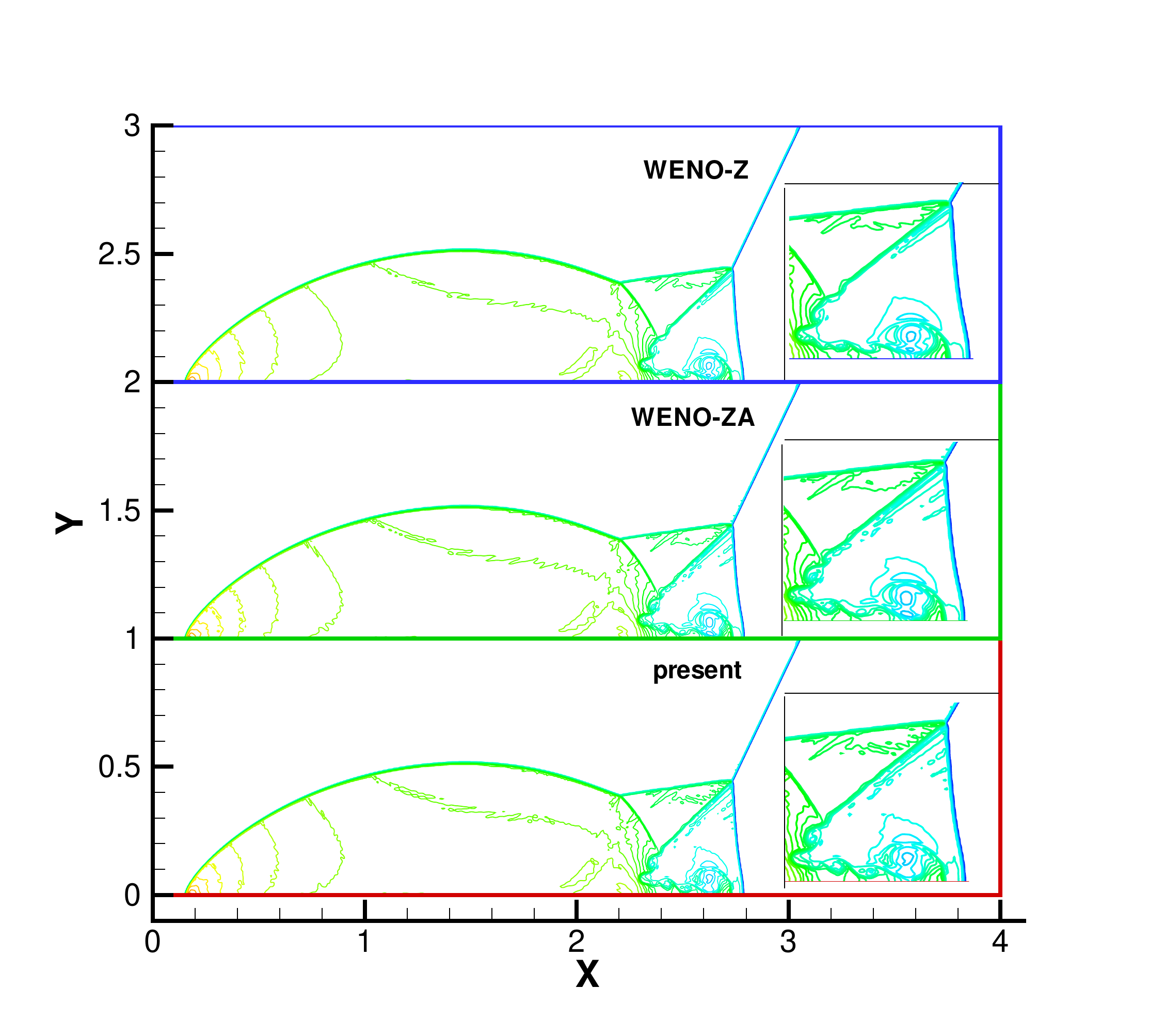}
\caption{Density contours of double Mach reflection, $960\times240$}
\label{fig:DM9}
\end{figure}

The computational domain is $[0, 4]\times [0, 1]$. For the bottom boundary, the exact post-shock condition is imposed for the
interval $[0, 0.6]$, and a reflective boundary condition is used for the rest. The top boundary is set to describe the exact
motion of a Mach 10 shock. Inflow and outflow boundary conditions are used for the left and right boundaries, respectively.
Fig.\ref{fig:DM9} gives the density contours on a mesh of $960\times 240$ at $t = 0.2$. The magnification of the roll-up region around the double Mach stems of each plot is also shown in the picture. It can be found that all schemes can capture shock structures well. However, from the three enlarged plots, we can see that the present scheme resolves
the roll-up structures more clearly than the other two schemes.

\subsubsection{Forward facing step flow}
This final test case is a two-dimensional flow past a forward facing step. It is usually used to show that increasing the resolution of a scheme can improve the
ability of capturing important details such as the roll-up of the vortex sheet via Kelvin-Helmholtz instability \cite{Cockburn1998, Balsara2009}. Our purpose is to prove that the new scheme performs robustly on this stringent
problem. Same as Ref.\cite{Balsara2009}, the two-dimensional wind tunnel spans a domain $[0,3]\times[0,1]$, and a forward facing step is set up at the coordinates $(0.6,0.2)$. The inflow boundary conditions are the ideal gas of Mach 3.0 with a density of 1.4 and a pressure of 1. The walls are set to be reflective boundaries. The ratio of specific
heats is 1.4.

\begin{figure}
\begin{minipage}{0.5\linewidth}
\centering
\includegraphics[width=1.0\textwidth]{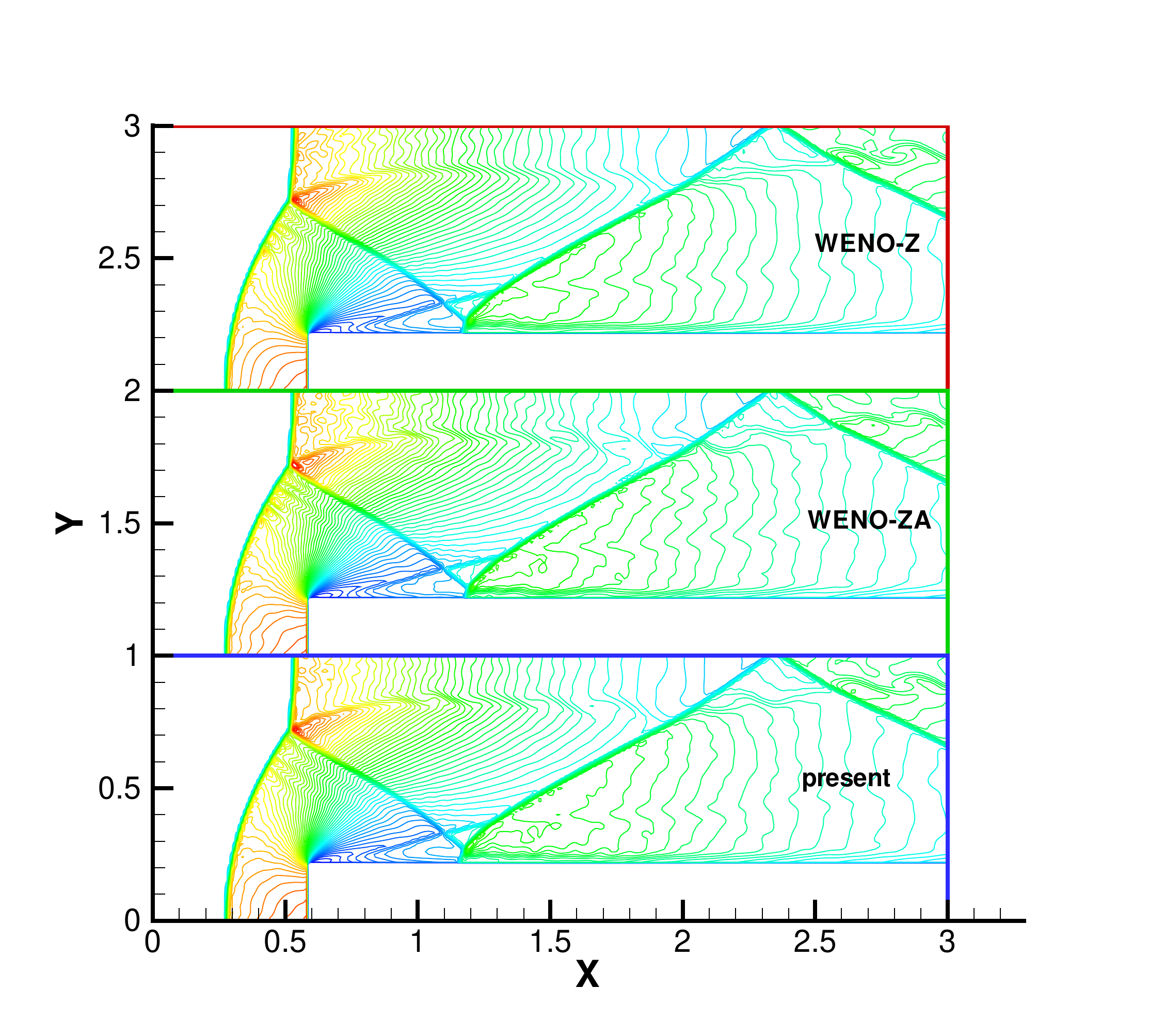}
\begin{minipage}{0.8\linewidth}
\caption{Density contours of forward facing step problem, $300\times 100$}
\label{fig:FS3}
\end{minipage}
\end{minipage}
\begin{minipage}{0.5\linewidth}
\centering
\includegraphics[width=1.0\textwidth]{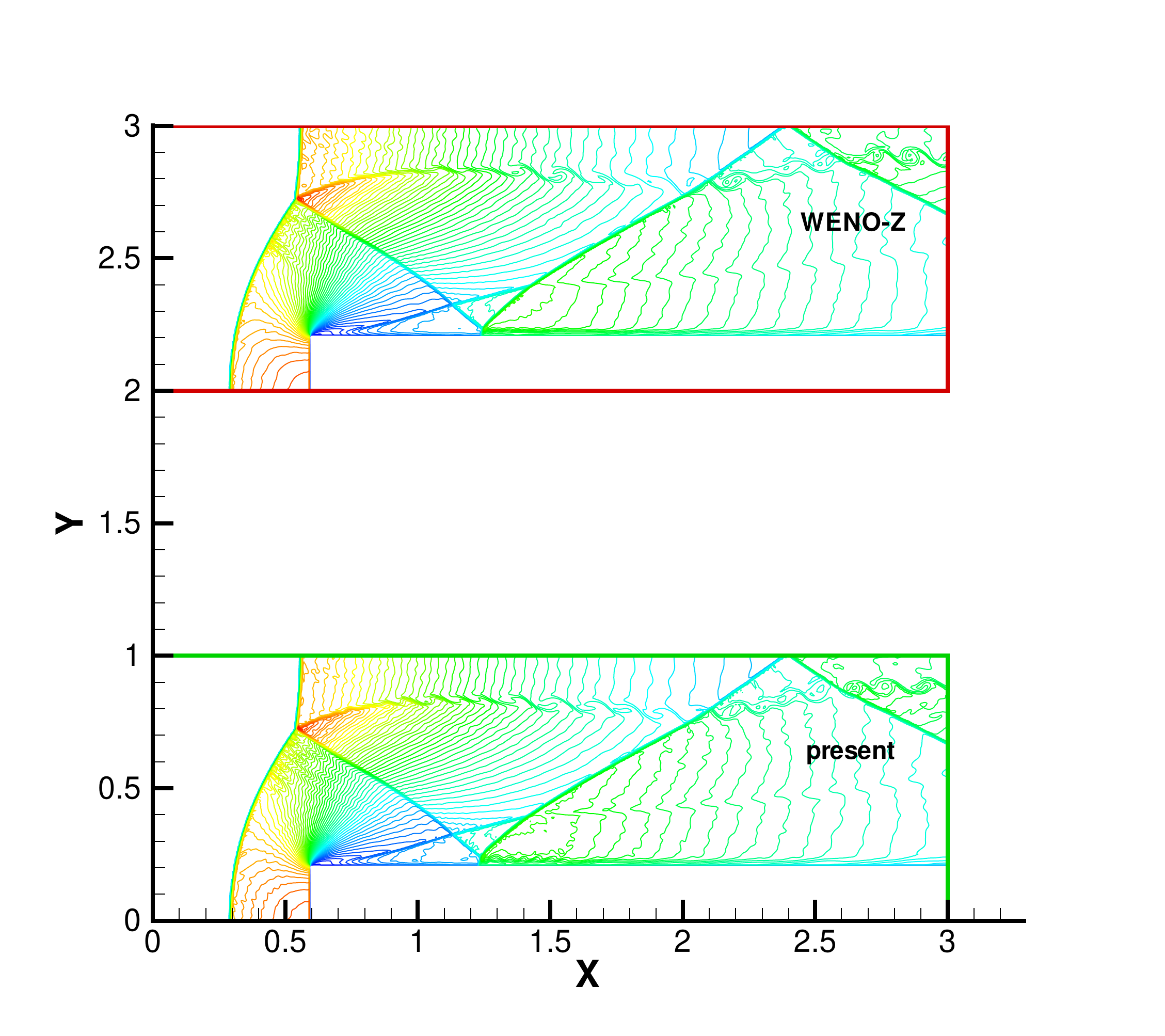}
\begin{minipage}{0.8\linewidth}
\caption{Density contours of forward facing step problem, $600\times 200$}
\label{fig:FS6}
\end{minipage}
\end{minipage}
\end{figure}

Two set meshes, i.e., $300\times 100$ and $600\times 200$, are tested. Figs.\ref{fig:FS3} and \ref{fig:FS6} are the density contours at the final time t=4.0. On the coarse $300\times 100$ mesh, all the shocks are properly captured, the WENO-ZA scheme and the present scheme obtain more clear instable structures than the WENO-Z scheme. On the finer mesh of $600\times 200$, the computation of WENO-ZA blows up and hence no result is obtained. From Fig.\ref{fig:FS6}, it can be seen that, the present scheme gives more clear roll-up of the vortex sheet. The computations of this case also show that, with the same computational conditions, the present scheme is more robust than the WENO-ZA scheme and more accurate and less dissipative than the WENO-Z scheme.

\section{Conclusion Remarks}
This paper presents a novel method for constructing WENO-Z type schemes. The method is mainly based on the analysis: in the formula for calculating the un-normalized weights of the fifth-order WENO-Z scheme, in order to capture shocks robustly, a relatively small value can be used to replace the constant $1$; on the contrary, in order to improve accuracy and reduce dissipation, it is beneficial to use a large value to replace $1$. Hence, first, we design a function of the local smoothness indicators of candidate sub-stencils to replace the constant $1$. The function can adaptively approach to a small value if the global stencil contains a discontinuity and approach to a large value if the global stencil is sufficiently smooth. Then, we suggest taking the square of the approximation of the fourth-order derivative, which is the maximal-order derivative can be approximated on a five-point stencil (the global stencil), as the global smoothness indicator.

Numerical results show that the new WENO-Z type scheme can achieve fifth-order accuracy at first-order critical point and fourth-order accuracy at second-order critical point. The new scheme has low numerical dissipation and is robust for solving problems with shocks.

This method can be easily extended to construct higher order WENO-Z type schemes which will be reported in an upcoming paper.

\section{Acknowledgement}
This research work was supported by the National Natural Science Foundation of China under Grants 11872067 and 91852203, NKRDPC 2016YFA0401200 and SCP No.TZ2016002.

\section*{References}
\bibliography{reference}
\end{document}